\documentclass[12pt]{article}
\usepackage{amssymb,amsmath,epsfig}

\newtheorem{theorem}{Theorem}
\newtheorem{proposition}[theorem]{Proposition}
\newtheorem{lemma}[theorem]{Lemma}

\newtheorem{corollary}[theorem]{Corollary}
\newenvironment{proof}{{\bf Proof }}{\par}

\newcommand{\Z}{{\mathbb Z}}
\newcommand{\Q}{{\mathbb Q}}
\newcommand{\R}{{\mathbb R}}
\newcommand{\C}{{\mathbb C}}
\newcommand{\CP}{{\mathbb P}}
\newcommand{\g}{{\mathfrak g}}
\newcommand{\h}{{\mathfrak h}}
\newcommand{\ct}{{\mathfrak t}}

\begin{document}

\title{Equivariant cohomology and equivariant intersection theory}

\author{Michel Brion}

\maketitle

{\sl Abstract.} This text is an introduction to equivariant
cohomology, a classical tool for topological transformation groups,
and to equivariant intersection theory, a much more recent topic
initiated by D. Edidin and W. Graham. It is based on lectures given at
Montr\'eal is Summer 1997.

Our main aim is to obtain explicit descriptions of cohomology or
Chow rings of certain manifolds with group actions which arise in
representation theory, e.g. homogeneous spaces and their
compactifications.

As another appplication of equivariant intersection theory, we
obtain simple versions of criteria for smoothness or rational
smoothness of Schubert varieties (due to Kumar \cite{Kumar},
Carrell-Peterson \cite{C} and Arabia \cite{Ar3}) whose statements and
proofs become quite transparent in this framework.
\bigskip

We now describe in more detail the contents of these notes; 
the prerequisites are notions of algebraic topology, compact Lie 
groups and linear algebraic groups.
Sections 1 and 2 are concerned with actions of compact Lie groups on
topological spaces, especially on symplectic manifolds. The material
of Section 1 is classical: universal bundles, equivariant cohomology
and its relation to usual cohomology, and the localization theorem for
actions of compact tori. 

A useful refinement of the latter theorem is presented in Section 2,
based on joint work with Mich\`ele Vergne. It leads to a complete
description of the cohomology ring of compact multiplicity-free
spaces. Examples include coadjoint orbits, projective toric manifolds
and De Concini-Procesi's complete symmetric varieties \cite{DP}.

The subject of the last three sections is equivariant intersection
theory for actions of linear algebraic groups on schemes. Edidin and
Graham's equivariant Chow groups are introduced in Section
3, after a brief discussion of usual Chow groups. The basic properties
of these equivariant groups are presented, as well as their
applications to Chow groups of quotients (after \cite{EG2}) and their
relation to usual Chow groups (after \cite{Br3}). As a direct
application, we obtain a structure theorem for the rational Chow ring
of a homogeneous space under a linear algebraic group; it turns out to
be much simpler than the corresponding rational cohomology ring.

Section 4, based on \cite{Br3}, presents Edidin and Graham's
localization theorem concerning equivariant Chow groups for torus
actions, and its refined version for projective smooth varieties. Then
we introduce  equivariant multiplicities at isolated fixed points, we
relate this notion to work of Joseph and Rossmann \cite{Joseph},
\cite{Ross} in representation theory, and we determine equivariant
multiplicities of Schubert varieties.

This is applied in Section 5 to criteria for smoothness or
rational smoothness of Schubert varieties, in terms of their equivariant
multiplicities. We deduce these criteria from new, more general
results concerning (rational) smoothness at an ``attractive'' fixed
point of a torus action; these results are closely related to recent
work of Arabia \cite{Ar3}.

The exposition follows the original five lectures, except in Section
5 which replaces a lecture on operators of divided differences, based
on \cite{Br3}. Indeed, we chose to develop here the applications of
equivariant intersection theory to singularities of Schubert
varieties, thereby answering questions raised by several participants
of the summer school. Further applications to spherical varieties will
be given elsewhere.

\section{Equivariant cohomology}

Let $X$ be a topological space and let $k$ be a commutative ring. We
denote by $H^*(X,k)$ the cohomology ring of $X$ with coefficients in
$k$. We will be mainly interested in the case where $k=\Q$ is the
field of rational numbers; we set $H^*(X):=H^*(X,\Q)$. 

Recall that $H^*(X,k)$ is a graded $k$-algebra, depending on $X$ in a
contravariant way: Any continuous map $f:X\to Y$ to a topological
space $Y$ induces a $k$-algebra homomorphism of degree zero
$f^*:H^*(Y,k)\to H^*(X,k)$ which only depends on the homotopy class of
$f$.

Consider now a topological space $X$ endowed with the action of a
topological group $G$; we will say that $X$ is a $G$-space. There
exists a principal $G$-bundle 
$$
p:E_G\to B_G
$$ 
such that the space $E_G$ is contractible, and such a bundle is
universal among principal $G$-bundles (see e.g. \cite{Hu} Chapter IV;
a construction of universal bundles will be given below when $G$ is a
compact Lie group). The group $G$ acts diagonally on $X\times E_G$ and
the quotient 
$$
X\times E_G\to (X\times E_G)/G=:X\times_G E_G
$$
exists. Define the {\sl equivariant cohomology ring} $H^*_G(X,k)$ by
$$
H^*_G(X,k):=H^*(X\times_G E_G,k).
$$
In particular, for $X$ a point, we have 
$$
H^*_G(pt,k)=H^*(E_G/G,k)=H^*(B_G,k).
$$
Because $p:E_G\to B_G$ is a principal bundle, the projection
$$
p_X:X\times_G E_G\to E_G/G=B_G
$$
is a fibration with fiber $X$. Thus, $H^*_G(X,k)$ is an algebra over
$H^*_G(pt,k)$. Moreover, restriction to the fiber of $p_X$ defines an
algebra homomorphism 
$$
\rho:H^*_G(X,k)/(H^+_G(pt,k))
=H^*_G(X,k)\otimes_{H^*_G(pt,k)}k\to H^*(X,k)
$$
where $(H^+_G(pt,k))$ denotes the ideal of $H^*_G(X,k)$ generated by
images of homogeneous elements of $H^*_G(pt,k)$ of positive degree.

We will see that, for certain spaces $X$ and for rational
coefficients, the map $\rho$ is an isomorphism, and that the 
equivariant cohomology algebra $H^*_G(X)$ can be described completely;
as a consequence, we will obtain a description of the usual cohomology
algebra of $X$. 

But it may happen that the map $\rho$ is trivial. Consider for example
a compact connected Lie group $G$ acting on $X=G$ by left
multiplication. Then 
$$
H^*_G(X)=H^*((G\times E_G)/G)=H^*(E_G)=\Q
$$
whereas $H^*(X)=H^*(G)$ is an exterior algebra on $l$ generators,
where $l$ is the rank of $G$ (a classical theorem of Hopf
\cite{Hopf}).

\newpage

\noindent
{\bf Remarks} 

1) The equivariant cohomology ring is independent of
the choice of a universal bundle. Indeed, if $E'_G\to B'_G$ is another
such bundle, then the projections $(X\times E_G\times E'_G)/G\to
X\times_G E_G$ and $(X\times E_G\times E'_G)/G\to X\times_G E'_G$ are
fibrations with contractible fibers $E'_G$ and $E_G$. Thus, both
projections induce isomorphisms in cohomology.

2) If $G$ acts trivially on $X$, then $X\times_G E_G=X\times B_G$.
By the K\"unneth isomorphism, it follows that
$$
H^*_G(X)\simeq H^*_G(pt)\otimes H^*(X).
$$
In particular, the $H^*_G(pt)$-module $H^*_G(X)$ is free.

On the other hand, if $G$ acts on $X$ with a quotient space
$X/G$ and with finite isotropy groups, then
$$
H^*_G(X)\simeq H^*(X/G).
$$
Indeed, the quotient map induces a map 
$$
\pi:X\times_G E_G\to X/G
$$
whose fibers are the quotients $E_G/G_x$ where $G_x$ denotes the
isotropy group of $x\in X$. Because $G_x$ is finite and $E_G$ is
contractible, the fibers of $\pi$ are $\Q$-acyclic. Thus, 
$\pi^*:H^*(X/G)\to H^*_G(X)$ is an isomorphism.

3) Let $H$ be a closed subgroup of $G$. Then the quotient $E_G/H$
exists and the map $E_G\to E_G/H$ is a universal bundle for $H$. As a
consequence, we obtain the following description of the
$G$-equivariant cohomology ring of the homogeneous space $G/H$:
$$
H^*_G(G/H,k)\simeq H^*_H(pt,k).
$$
Indeed, we have
$$H^*_G(G/H,k)=H^*(G/H\times_G E_G,k)=H^*(E_G/H,k).
$$
More generally, let $Y$ be a $H$-space. The quotient of $G\times Y$ by
the diagonal action of $H$ exists; we denote it by $G\times_H Y$. Then
we obtain as above:
$$
H^*_G(G\times_H Y,k)\simeq H^*_H(Y,k).
$$

4) Let $n$ be a positive integer. Let $E_G^n\to B_G^n$ be a principal
$G$-bundle such that $E_G^n$ is $n$-connected (that is, any continuous
map from a sphere of dimension at most $n$ to $E_G^n$ is homotopic to
the constant map). Then 
$$
H^m_G(X,k)=H^m(X\times_G E_G^n,k)
$$
for any compact topological $G$-manifold $X$ of dimension at most $n$,
and for any integer $m\leq n$ \cite{Hu} Chapter IV, Theorem 13.1. In
other words, we can compute any homogeneous component of $H^*_G(X)$ by
finite approximations of $E_G$. This will be the starting point for
the definition of equivariant Chow groups, in Section 3.
\bigskip

\noindent
{\bf Examples} 

1) $G=S^1$ (the multiplicative group of complex numbers
of modulus one). Let $G$ act on $\C^n$ by scalar multiplication; then
the unit sphere $S^{2n-1}$ is $G$-stable and the quotient $S^{2n-1}/G$
is complex projective space $\CP^{n-1}$. We thus obtain a principal
bundle $S^{2n-1}\to\CP^{n-1}$ with total space $(2n-2)$-connected. So
we can take $E_G=\lim_{\to}S^{2n-1}=S^{\infty}$; then
$B_G=\lim_{\to}\CP^n=\CP^{\infty}$. It follows that 
$$
H^*_G(pt,k)=H^*(B_G,k)=k[t]
$$
where $t$ is an indeterminate of degree 2.

2) More generally, let $G$ be a (compact) torus. Then 
$G\simeq (S^1)^{\ell}$ and thus, we can take $E_G=(E_{S^1})^{\ell}$.
In particular, $H^*_G(pt,k)$ is a polynomial $k$-algebra on $\ell$
indeterminates of degree 2. A more intrinsic description of
$H^*_G(pt)$ is as follows. Denote by $\Xi(G)$ the character group of
$G$ consisting of all continuous group homomorphisms $G\to S^1$. Any
$\chi\in\Xi(G)$ defines a one-dimensional complex representation of
$G$ with space $\C\chi$. Consider the associated complex line bundle
$$
L(\chi):=(E_G\times_G \C\chi \to B_G)
$$
and its first Chern class $c(\chi)\in H^2(B_G)$. Let $S$ be the
symmetric algebra over $\Q$ of the group $\Xi(G)$. Then $S$ is
a polynomial ring on $\ell$ generators of degree 1, and the map 
$\chi\to c(\chi)$ extends to a ring isomorphism
$$
c:S\to H^*_G(pt)
$$
which doubles degrees: the {\sl characteristic homomorphism}. 

3) Finally, let $G$ be a compact Lie group. Then we can embed $G$ as a
closed subgroup of some $GL_n$, where $GL_n$ denotes the group of
invertible $n\times n$ {\sl complex} matrices. We now construct a
universal bundle $E_{GL_n}\to B_{GL_n}$ for the group $GL_n$; then
$E_{GL_n}\to E_{GL_n}/G$ will be a universal bundle for $G$.

Let $N\geq n$ be an integer. Consider the space $M_{N\times n}$ of all
$N\times n$ complex matrices, and its subset $M_{N\times n}^{max}$ of
matrices of maximal rank $n$. Clearly, $GL_n$ acts on $M_{N\times n}$
and preserves $M_{N\times n}^{max}$. Let $Gr_{N,n}$ be the Grassmann
variety of linear subspaces of dimension $n$ of $\C^N$. With any
matrix $A\in M_{N\times n}^{max}$, we associate its image $Im(A)$.
This defines a map
$$
M_{N\times n}^{max}\to Gr_{N,n}
$$
which is a principal $GL_n$-bundle. Moreover, the complement of
$M_{N\times n}^{max}$ in $M_{N\times n}$ is a (Zariski) closed
subset of $M_{N\times n}$ of codimension $N-n+1$ over $\C$. It follows
that $M_{N\times n}^{max}$ is $(N-n)$-connected. So we obtain a
construction of $E_{GL_n}$ as the space of $\infty\times n$ complex
matrices of maximal rank.

\bigskip

We now show that equivariant cohomology for a compact connected Lie
group can be described in terms of equivariant cohomology for a
maximal torus.

\begin{proposition}\label{deg}
Let $G$ be a compact connected Lie group and let $T\subset G$ be a
maximal torus with normalizer $N$ and with Weyl group $W=N/T$; let $X$
be a $G$-space.

(i) The group $W$ acts on $H^*_T(X)$ and we have an isomorphism
$$
H^*_G(X)\simeq H^*_T(X)^W.
$$
In particular, $H^*_G(pt)$ is isomorphic to $S^W$ where $S$ denotes
the symmetric algebra of the character group $\Xi(T)$ (occuring in
degree 2), and $S^W$ the ring of $W$-invariants in $S$.

(ii) The map
$$
S\simeq H^*_G(G/T)\to H^*(G/T)
$$
is surjective and induces an isomorphism $S/(S^W_+)\to H^*(G/T)$
where $(S^W_+)$ denotes the ideal of $S$ generated by all homogeneous
$W$-invariants of positive degree.

(iii) We have an isomorphism
$$
S\otimes_{S^W} H^*_G(X)\simeq H^*_T(X).
$$
In particular, $H^*_T(G/T)$ is isomorphic to $S\otimes_{S^W}S$.
\end{proposition}
\begin{proof} (after \cite{Hsiang} Chapter III, \S 1).

(i) Let $G^{\C}$ be the complexification of $G$; then $G^{\C}$ is a
complex connected reductive Lie group. Let $B$ be a Borel subgroup of
$G^{\C}$ containing the compact torus $T$. Then, by the Iwasawa
decomposition, we have $G^{\C}=GB$ and $G\cap B=T$. Thus, the map
$G/T\to G^{\C}/B$ is a homeomorphism. By the Bruhat decomposition, the
flag manifold $G^{\C}/B$ has a stratification by $\vert W\vert$ strata,
each of them being isomorphic to a complex affine space. It follows
that $H^*(G/T)$ vanishes in odd degrees, and that the topological
Euler characteristic $\chi(G/T)$ is equal to $\vert W\vert$. 

Because the finite group $W$ acts freely on $G/T$ with quotient $G/N$,
we have an isomorphism
$$
H^*(G/N)\simeq H^*(G/T)^W
$$
and moreover, $\chi(G/N)=\vert W\vert^{-1}\chi(G/T)=1$. It follows
that $H^*(G/N)$ vanishes in odd degrees and is one-dimensional. In
other words, $G/N$ is $\Q$-acyclic. 

Thus, the fibration 
$$
X\times_N E_G\to X\times_G E_G
$$
with fiber $G/N$ induces an isomorphism $H^*_G(X)\to H^*_N(X)$.
Moreover, we have a covering 
$$
X\times_T E_G\to X\times_N E_G
$$
with group $W$, whence $H^*_N(X)\simeq H^*_T(X)^W$.

(ii) By (i) applied to the point, the odd degree part of
$$
H^*(B_G)=H^*_G(pt)=S^W
$$ 
is zero. The fibration 
$$
G/T\times_G E_G\to B_G
$$
has fiber $G/T$ with vanishing odd cohomology, too. Thus, the Leray
spectral sequence degenerates, and we obtain an isomorphism
$$
H^*(G/T)\simeq H^*_G(G/T)/(H^*(B_G)_+)=H^*_T(pt)/(H^*_T(pt)^W_+)
=S/(S^W_+).
$$

(iii)
For the fibration 
$$
X\times _T E_G\to X\times _G E_G
$$
with fiber $G/T$, restriction to fiber defines a ring homomorphism
$H^*_T(X)\to H^*(G/T)$ which is surjective, because the map 
$S\to H^*(G/T)$ is. Thus, the Leray spectral sequence degenerates, and
$H^*_T(X)\simeq H^*_G(X)\otimes H^*(G/T)$ as a $H^*_G(X)$-module. This
implies our statement.
\end{proof}

\bigskip

\noindent
{\bf Remark.} The map $S\to H^*(G/T)$ is the {\sl characteristic
homomorphism}, which associates with any $\chi\in\Xi(T)$ the
Chern class of the corresponding line bundle on $G/T$. The map
$S\otimes_{S^W}S\to H^*_T(G/T)$ is a lift of this characteristic
homomorphism to $T$-equivariant cohomology.
\bigskip

For $G$ and $X$ as before, consider the restriction map 
$$
\rho:H^*_G(X)/(S^W_+)\to H^*(X).
$$ 
We already saw that $\rho$ may not be surjective; however, there is an
important class of $G$-spaces for which $\rho$ is an isomorphism.
\bigskip

\noindent
{\bf Definition.} Let $G$ be a compact connected Lie group and let $X$
be a $G$-space. Then $X$ is a {\sl symplectic $G$-manifold} if $X$ is
a smooth manifold with a non-degenerate closed, $G$-invariant form
$\omega$ (the symplectic form). Moreover, $X$ is a {\sl Hamiltonian
$G$-space} if it satisfies the following additional condition:

There exists a $G$-equivariant smooth map $\mu:X\to\g^*$ (where
$\g$ denotes the Lie algebra of $G$) such that, for any $x\in X$,
$\xi\in T_x X$ and $\eta\in\g$, we have: 
$$
d\mu_x(\xi)(\eta)=\omega_x(\xi,\eta_x)
$$
(where $\eta_x$ is the value at $x$ of the vector field on $X$ defined
by $\eta\in\g$).
\bigskip

The map $\mu$ is called a {\sl moment map} for the action of
$G$ on $X$; it is unique up to translation by a central element of
$\g$. A Hamiltonian $G$-space $X$ is also Hamiltonian for the action
of any closed connected subgroup $H$ of $G$; the moment map for this
action is $\mu$ followed by restriction $\g^*\to\h^*$. 

Examples of compact Hamiltonian $G$-spaces are complex projective
manifolds with a linear action of $G$. More precisely, let $V$ be the
space of a finite-dimensional complex representation of $G$, let
$\CP(V)$ be its projectivization, and let $X\subset\CP(V)$ be a
$G$-invariant projective submanifold. The choice of a $G$-invariant
Hermitian scalar product $(\cdot,\cdot)$ on $V$ defines a symplectic 
form on $\CP(V)$ (the imaginary part of the Fubini-Study metric). 
This makes $\CP(V)$ a Hamiltonian $G$-space with moment map given by
$$
\mu([v])(\eta)=\frac{(v,\eta v)}{2i\pi(v,v)}
$$
where $[v]\in\CP(V)$ is the image of $v\in V$, and $\eta\in\g$.
Moreover, $X$ is a Hamiltonian subspace of $\CP(V)$ \cite{Kirwan}.

Now we can state the following result, due to F. Kirwan \cite{Kirwan}
Proposition 5.8.

\begin{proposition}\label{kir}
Let $G$ be a compact connected Lie group and let $X$ be a compact
Hamiltonian $G$-space. Then, notation being as above, the $S^W$-module
$H^*_G(X)$ is free, and the map $\rho:H^*_G(X)/(S^W_+)\to H^*(X)$ is
an isomorphism.
\end{proposition}

As an example, consider the $G$-space $X=G/T$ where, as before, $T$ is
a maximal torus of $G$. Embed $G/T$ into $\g^*$ as the orbit of a
regular element. Then, as a coadjoint orbit, $G/T$ is a Hamiltonian
$G$-manifold and the moment map is simply inclusion $G/T\to\g^*$, see
e.g. \cite{Au} Chapter II, 3.3. In this way, we recover the structure
of the cohomology ring of the flag manifold $G/T$.
\bigskip

A powerful tool in equivariant cohomology is the following
localization theorem of Borel-Atiyah-Segal \cite{Hsiang} Chapter III,
\S 2.

\begin{theorem}\label{loc}
Let $T$ be a compact torus and let $X$ be a $T$-space which embeds
equivariantly into a finite-dimensional $T$-module. Let 
$i_T:X^T\to X$ be the inclusion of the fixed point set. Then the 
$S$-linear map 
$$
i_T^*:H^*_T(X)\to H^*_T(X^T)
$$
becomes an isomorphism after inverting finitely many non-trivial
characters.
\end{theorem}
\begin{proof}
For later use, we prove the following stronger statement.
\begin{lemma}\label{loca}
Let $\Gamma\subset T$ be a closed subgroup and let 
$i_{\Gamma}:X^{\Gamma}\to X$ be the inclusion of the fixed point set.
If $X$ embeds equivariantly into a finite-dimensional $T$-module, then
the restriction map
$$
i_{\Gamma}^*:H^*_T(X)\to H^*_T(X^{\Gamma})
$$
becomes an isomorphism after inverting finitely many characters of $T$
which restrict non-trivially to $\Gamma$.
\end{lemma}
\begin{proof} 
First we consider the case where $X^{\Gamma}$ is empty. Then we have
to prove that $H^*_T(X)$ is killed by a product of characters of $T$
which restrict non-trivially to $\Gamma$. Because
$X$ embeds into a finite-dimensional $T$-module, we can cover $X$
by finitely many $T$-invariant subsets $X_1,\ldots,X_n$ such that
every $X_j$ admits a $T$-equivariant map $X_j\to T/\Gamma_j$ where
$\Gamma_j\subset T$ is a closed subgroup which does not contain
$\Gamma$ (exercise). Thus, we can choose $\chi_j\in\Xi(T)$ such that
$\chi_j$ restricts trivially to $\Gamma_j$ but not to $\Gamma$. For
any $T$-invariant subset $U_j\subset X_j$, the ring $H^*_T(U_j)$ is a
module over $H^*_T(T/\Gamma_j)$. The image of $\chi_j$ in the latter
is zero; thus, multiplication by $\chi_j$ is zero in $H^*_T(U_j)$. By
Mayer-Vietoris, multiplication by the product of the $\chi_j$ is zero
in $H^*_T(X)$.

In the general case, let $Y\subset X$ be a closed $T$-stable
neighborhood of $X^{\Gamma}$ in $X$. Let $Z$ be the closure of
$X\setminus Y$ in $X$; then $Z$ is $T$-stable and $Z^{\Gamma}$ is
empty. By the first step of the proof and Mayer-Vietoris, it follows
that restriction $H^*_T(X)\to H^*_T(Y)$ is an isomorphism after
inverting a finite set ${\cal F}$ of characters of $T$. Moreover,
${\cal F}$ is independent of $Y$, because it only depends on a
$T$-module which contains $X$. To conclude the proof,
observe that $H^*_T(X^{\Gamma})$ is the direct limit of the
$H^*_T(Y)$.
\end{proof}
\end{proof}
\bigskip

By \cite{Kawa} Theorem 4.12, the assumptions of Theorem \ref{loc} hold
if $X$ is a compact manifold. If moreover $X$ is a symplectic
$T$-manifold, then $X^T$ is, too, as a consequence of the following 

\begin{lemma}\label{fix}
Let $G$ be a compact connected Lie group, $\Gamma\subset G$ a closed
subgroup with centralizer $G^{\Gamma}$, and $X$ a symplectic
$G$-manifold. Then the fixed point set $X^{\Gamma}$ is a symplectic
$G^{\Gamma}$-manifold, and the normal bundle $N_{X,X^{\Gamma}}$ has a
natural structure of a complex vector bundle. 

If moreover the $G$-action on $X$ is Hamiltonian with moment map
$\mu$, then the $G^{\Gamma}$-action on $X^{\Gamma}$ is Hamiltonian
with moment map: restriction of $\mu$ followed by restriction to
$\g^{\Gamma}$.
\end{lemma}
\begin{proof}
Let $x\in X^{\Gamma}$. Because $\Gamma$ is compact, there exists a
$\Gamma$-invariant neighborhood of $x$ in $X$ which is equivariantly
isomorphic to a $\Gamma$-module \cite{Kawa} Theorem 4.10. It follows
that $X^{\Gamma}$ is a closed submanifold of $X$, and that we have the
equality of tangent spaces
$$
T_x(X^{\Gamma})=(T_x X)^{\Gamma}.
$$
Moreover, $(T_x X)^{\Gamma}$ has a unique $\Gamma$-invariant
complement in $T_x X$: the span of all $\gamma\xi-\xi$ where
$\gamma\in\Gamma$ and $\xi\in T_x X$. Because the symplectic form
$\omega_x$ is $\Gamma$-invariant, it follows that $(T_x X)^{\Gamma}$
is a symplectic subspace of $T_x X$. Thus, restriction
$\omega^{\Gamma}$ of $\omega$ to $X^{\Gamma}$ is non-degenerate, and
closed because $\omega$ is. 

If moreover $X$ is Hamiltonian with moment map $\mu$, then we have for
$\xi\in(T_x X)^{\Gamma}$ and $\eta\in\g^{\Gamma}$:
$$
\omega^{\Gamma}_x(\xi,\eta_x)=\omega_x(\xi,\eta_x)=d\mu_x(\xi)(\eta),
$$
that is, $\mu\vert_{X^{\Gamma}}$ is a moment map for the action of
$G^{\Gamma}$.
\end{proof}
\bigskip

For a symplectic $T$-manifold $X$, this leads to an explicit form
of the localization theorem. Let $Y$ be a connected component of
$X^T$, then the codimension of $Y$ in $X$ is an even number, say $2n$,
and we have a Gysin map
$$
i_{Y*}: H^*_T(Y)\to H^*_T(X)
$$
which increases degree by $2n$. Composition $i_Y^*\circ i_{Y*}$ is
multiplication by the equivariant Euler class of $Y$ in $X$, which we
denote by $Eu_T(Y,X)$: because the normal bundle $N_{Y,X}$ is a
complex equivariant vector bundle, it admits equivariant Chern classes
in $H^*_T(Y)$ (which are the Chern classes of the vector bundle
$(Y\times N_{Y,X}\times E_T)/T$ over $Y_T$) and $Eu_T(Y,X)$ is the top
equivariant Chern class
$c_{n,T}(N_{X,Y})$. Moreover, each fiber of $N_{X,Y}$ is a complex
$T$-module with non-zero weights. Denoting by $det(N_{Y,X})$ the
product of these weights, we have in $H^*_T(Y)=S\otimes H^*(Y)$:
$$
c_{n,T}(N_{Y,X})\equiv det(N_{Y,X})\otimes[Y]
\eqno(mod~S\otimes H^+(Y)).
$$
Because $S\otimes H^+(Y)$ consists of nilpotent elements,
$c_{n,T}(N_{Y,X})$ becomes invertible in $H^*_T(Y)$ after inverting
$det(N_{Y,X})$. Moreover, setting $i_{T*}:=\sum_Y i_{Y*}$, we have
$$
i_T^*(i_{T*}(u))=\sum_Y c_{top,T}(N_{Y,X})\cup u
$$
for any $u\in H^*_T(X^T)$. As a consequence, $i_{T*}$ is an
isomorphism after localization, with inverse 
$$
u\mapsto\sum_Y c_{top,T}(N_{Y,X})^{-1}\cup i_Y^*(u).
$$

This ``covariant'' version of the localization theorem will be
extended in Section 4 to equivariant Chow groups, which can be
considered as algebraic equivariant {\sl homology} groups. It has
many important applications to residue or integral formulae, in
particular the Bott residue formula and the Duistermaat-Heckman
theorem, for which we refer to \cite{AB} and \cite{Au}.

\section{A precise version of the localization theorem}

As before, we consider a compact torus $T$ and a $T$-space $X$; we
denote by $i_T:X^T\to X$ the inclusion of the fixed point set. We wish
to describe the image of the restriction map 
$$
i_T^*:H^*_T(X)\to H^*_T(X^T).
$$
For this, let $T'\subset T$ be a subtorus; observe that
$$
i_T:X^T\to X
$$
factors as 
$$
i_{T,T'}:X^T\to X^{T'}
$$
followed by 
$$
i_{T'}:X^{T'}\to X.
$$
Thus, the image of $i_T^*$ is contained in the image of $i_{T,T'}^*$.
This observation will lead to a complete description of the image of
$i_T^*$, if the $S$-module $H^*_T(X)$ is free. Indeed, we have the
following result, due to Chang-Skjelbred \cite{GS} and Hsiang
\cite{Hsiang} Corollary p. 63, in a slightly different formulation.

\begin{theorem}\label{locprec}
Let $X$ be a $T$-space which admits an equivariant embedding into the
space of a finite-dimensional representation of $T$. If the $S$-module
$H^*_T(X)$ is free, then the map
$$
i_T^*:H^*_T(X)\to H^*_T(X^T)
$$
is injective, and its image is the intersection of the images of the
maps 
$$
i_{T,T'}^*:H^*_T(X^{T'})\to H^*_T(X^T)
$$
where $T'$ runs over all subtori of codimension 1 of $T$.
\end{theorem}
\begin{proof}
We already know that $i_T^*$ becomes an isomorphism after inverting a
finite family ${\cal F}$ of non-trivial characters of $T$. Because the
$S$-module $H^*_T(X)$ is free, it follows that $i_T^*$ is injective.

It remains to prove that the intersection of the images of the
$i_{T,T'}^*$ is contained in the image of $i_T^*$. Choose a basis 
$(e_j)_{j\in J}$ of the free $S$-module $H^*_T(X)$.
For any $j\in J$, let 
$$
e_j^*:H^*_T(X)\to S
$$
be the corresponding coordinate function. Then there exists a
$S$-linear map 
$$
f_j:H^*_T(X^T)\to S[1/\chi]_{\chi\in{\cal F}}
$$
such that $f_j\circ i_T^*=e_j^*$.

We may assume that each $\chi\in{\cal F}$ is primitive, i.e., not
divisible in $\Xi(T)$. Then its kernel 
$ker(\chi)\subset T$ is a subtorus of codimension 1. Let $u$ be in the
image of $i_{T,ker(\chi)}^*$; write 
$$
u=i_{T,ker(\chi)}^*(v)
$$ 
where $v\in H^*_T(X^{ker(\chi)})$. By Lemma \ref{loca} applied to 
$\Gamma=ker(\chi)$, there exist a product $P_{\chi}$ of weights of $T$
which are not multiples of $\chi$, such that $P_{\chi} v$ is in the
image of $i^*_{ker(\chi)}$. It follows that $P_{\chi}u$ is in the
image of $i^*_T$. Applying $f_j$, we obtain $P_{\chi} f_j(u)\in S$.
Thus, the denominator of $f_j(u)$ is not divisible by $\chi$.

If $u\in H^*_T(X^T)$ is in the intersection of the images of the
$i^*_{T,ker(\chi)}$ for all $\chi\in{\cal F}$, then
$f_j(u)\in S[1/\chi]_{\chi\in{\cal F}}$ but the denominator of 
$f_j(u)$ is not divisible by any element of ${\cal F}$; whence
$f_j(u)\in S$. It follows that
$u=i^*_T(\sum_{j\in J}\,f_j(u) e_j)$ 
is in the image of $i^*_T$.
\end{proof}

\bigskip

Observe that the assumptions of Theorem \ref{locprec} are satisfied
whenever $X$ is a compact Hamiltonian $T$-space. Indeed, the
$S$-module $H^*_T(X)$ is free by Proposition \ref{kir}, and the fixed
point set of any closed subgroup of $T$ is a Hamiltonian $T$-space by
Lemma \ref{fix}. 

If moreover $X^T$ is finite and each $X^{T'}$ has dimension at most 2,
then each connected component $Y$ of $X^{T'}$ is either a point or a
2-sphere, see e.g. \cite{Au} Chapter I, 3.3. In the latter case, we
can see $Y$ as complex projective line where $T$ acts through
multiplication by a character $\chi$; then $Y^T=\{0,\infty\}$. It is
easily checked that
$$
i^*_TH^*_T(Y)=\{(f_0,f_{\infty})\in S\times S~\vert~
f_0\equiv f_{\infty}~(mod~\chi)\}.
$$

From the discussion above, we deduce the following explicit
description of the image of $H^*_T(X)$ under restriction to the fixed 
point set, for certain Hamiltonian $T$-spaces $X$ \cite{GKM}.

\begin{corollary}
Let $X$ be a compact Hamiltonian $T$-space with finitely many fixed
points $x_1,\ldots,x_m$ and such that $\dim(X^{T'})\leq 2$ for any
subtorus $T'\subset T$ of codimension 1. Then, via $i_T^*$, the
algebra $H^*_T(X)$ is isomorphic to the subalgebra of $S^m$ consisting
of all $m$-tuples $(f_1,\ldots,f_m)$ such that: 
$$
f_j\equiv f_k \eqno(mod~\chi)
$$
whenever the fixed points $x_j$ and $x_k$ are in the same connected
component of $X^{ker(\chi)}$, for $\chi$ a primitive character of $T$.
Moreover, the cohomology algebra $H^*(X)$ is the quotient of
$H^*_T(X)$ by its ideal generated by all $(f,f,\ldots,f)$ where 
$f\in S$ is homogeneous of positive degree. 
\end{corollary}

This statement can be reformulated in a more geometric way, by
introducing the affine algebraic scheme $V(X)$ over $\Q$ associated
with the finitely generated $\Q$-algebra $H^*_T(X)$. Indeed, $V(X)$ is
obtained from the disjoint union of $m$ copies of $\ct$ by identifying
the $j$-th and $k$-th copies along their common hyperplane $(\chi=0)$
whenever $x_j$ and $x_k$ are in the same connected component of
$X^{ker(\chi)}$. In particular, $V(X)$ is reduced. The inclusion
$S\subset H^*_T(X)$ defines a morphism $V(X)\to\ct$ which restricts
to the identity on each copy of $\ct$. This morphism is finite and
flat (because $H^*_T(X)$ is a finite free $S$-module) and its
scheme-theoretic fiber at the origin is the scheme associated with
$H^*(X)$.
\bigskip

\noindent
{\bf Examples} 

1) (flag manifolds) Let $G$ be a compact connected Lie
group, $T\subset G$ a maximal torus, and $X=G/T$ the flag manifold of
the complexification of $G$. Denote by $x$ the base point of $X$. Then
the fixed point set $X^T$ is the orbit $Wx$; it identifies to the Weyl
group $W$. Denote by $\Phi$ the root system of $(G,T)$. Let $\chi$ be
a primitive character of $T$. If $\chi$ is not in $\Phi$, then
$X^{ker(\chi)}=X^T$. If $\chi$ is in $\Phi$, with corresponding
reflection $s\in W$, then  
$$
X^{ker(\chi)}=G^{ker(\chi)}Wx
$$ 
is a disjoint copy of  $\vert W\vert/2$ complex projective lines
joigning the fixed points $wx$ and $swx$ for all $w\in W$. So we
obtain from the corollary:
$$
H^*_T(G/T)=\{(f_w)_{w\in W}~\vert~f_w\in S,~f_w\equiv f_{s_{\alpha}w}
~(mod~\alpha)~\forall \alpha\in \Phi, \forall w\in W\}.
$$
In other words, the scheme $V(G/T)$ is the union of $\vert W\vert$
copies of $\ct$ identified along all reflection hyperplanes.

On the other hand, we saw that $H^*_T(G/T)=S\otimes_{S^W} S$. This
can be related to the description above, as follows. The
homomorphism $S\otimes 1\to H^*_T(G/T)$ is given by the structure of
$S$-module of $H^*_T(G/T)$; its composition with $i_T^*$ maps each
$f\in S$ to the tuple $(f,f,\ldots,f)$. The homomorphism 
$1\otimes S\to H^*_T(G/T)$ is the characteristic homomorphism; its
composition with $i_T^*$ maps $f$ to the tuple $(w(f))_{w\in W}$. The
scheme $V(G/T)$ identifies to the union of subspaces of $\ct\times \ct$
which are images of the diagonal $\Delta(\ct)$ by some element of
$W\times W$. Indeed, we have
$$
V(X)=Spec(S\otimes_{S^W}S)=\ct\times_{\ct/W}\ct=
\{(a,b)\in\ct\times\ct~\vert~a\in Wb\}.
$$
This is also the union of $\vert W\vert$ copies of $\ct$ identified
along reflection hyperplanes.

For other pictures of the equivariant cohomology ring of the flag
manifold, we refer to \cite{Ar1}, \cite{Ar2} and \cite{KK}; see also
Section 4 and Peterson's lecture notes in this volume.

2) (toric manifolds) Let $X$ be a projective toric manifold, that is,
$X$ is a complex algebraic projective manifold and the
complexification $T^{\C}$ of $T$ acts on $X$ with a dense orbit
isomorphic to $T^{\C}$. Then $X$ is a Hamiltonian $T$-space,
$\mu\vert_{X^T}$ is injective, and $\mu(X):=P$ is a convex polytope in
$\ct^*$ with vertex set $\mu(X^T):={\cal V}$. Moreover, $\mu$ induces 
a bijection bewteen $T^{\C}$-orbit closures in $X$ and faces of $P$
\cite{Au} Chapter 4. It follows that, for a character $\chi$ of $T$,
the set $X^{ker(\chi)}$ is equal to $X^T$, except when there exists an
edge of $P$ with direction $\chi$. In this case, $X^{ker(\chi)}$ is a
complex projective line joigning the $T$-fixed points associated with
the vertices of this edge. Thus, $i_T^*H^*_T(X)$ consists of all
families $(f_v)_{v\in{\cal V}}$ of $S$, such that $f_v\equiv f_{v'}$
(mod $\chi$) whenever the segment $[v,v']$ is an edge of $P$ with
direction $\chi$.

Let $F$ be a face of $P$. Choose a point $x_F$ in the relative
interior of $F$, then the union of all half-lines with origin 
at $x_F$ and which meet $P\setminus x_F$ is a closed convex cone which
depends only on $F$; we denote this cone by $C_F P$. The translated
cone $-x_F+C_F P$ only depends on $F$, too, and contains the origin of
$\ct^*$; denote by $\sigma_F\subset \ct$ its dual cone. Then the cones
$\sigma_F$ are a subdivision of $\ct$, with maximal cones $\sigma_v$
($v\in{\cal V}$). 

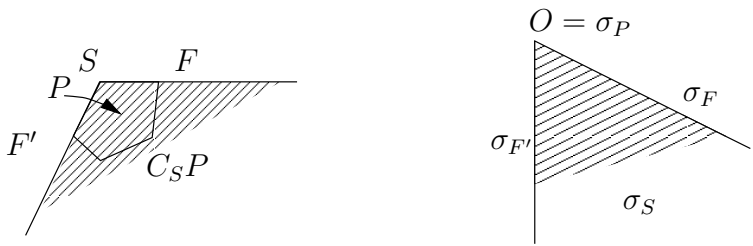
\begin{center}
\begin{picture}(0,0)%
\epsfig{file=figure1.pstex}%
\end{picture}%
\setlength{\unitlength}{0.00083300in}%
\begingroup\makeatletter\ifx\SetFigFont\undefined%
\gdef\SetFigFont#1#2#3#4#5{%
  \reset@font\fontsize{#1}{#2pt}%
  \fontfamily{#3}\fontseries{#4}\fontshape{#5}%
  \selectfont}%
\fi\endgroup%
\begin{picture}(4689,1503)(601,-936)
\put(4824,-24){\makebox(0,0)[lb]{\smash{\SetFigFont{12}{14.4}{\rmdefault}{\mddefault}{\updefault}$\sigma_{F}$}}}
\put(4449,-699){\makebox(0,0)[lb]{\smash{\SetFigFont{12}{14.4}{\rmdefault}{\mddefault}{\updefault}$\sigma_{S}$}}}
\put(3624,-324){\makebox(0,0)[lb]{\smash{\SetFigFont{12}{14.4}{\rmdefault}{\mddefault}{\updefault}$\sigma_{F'}$}}}
\put(1651,164){\makebox(0,0)[lb]{\smash{\SetFigFont{12}{14.4}{\rmdefault}{\mddefault}{\updefault}$F$}}}
\put(601,-361){\makebox(0,0)[lb]{\smash{\SetFigFont{12}{14.4}{\rmdefault}{\mddefault}{\updefault}$F'$}}}
\put(1051,164){\makebox(0,0)[lb]{\smash{\SetFigFont{12}{14.4}{\rmdefault}{\mddefault}{\updefault}$S$}}}
\put(856, -6){\makebox(0,0)[lb]{\smash{\SetFigFont{12}{14.4}{\rmdefault}{\mddefault}{\updefault}$P$}}}
\put(1476,-491){\makebox(0,0)[lb]{\smash{\SetFigFont{12}{14.4}{\rmdefault}{\mddefault}{\updefault}$C_SP$}}}
\put(3864,411){\makebox(0,0)[lb]{\smash{\SetFigFont{12}{14.4}{\rmdefault}{\mddefault}{\updefault}$O=\sigma_P$}}}
\end{picture}

\end{center}

We can consider each $f_v\in S$ as a polynomial function on
$\sigma_v$. Then the congruences above mean that the functions $f_v$
are compatible along all common facets of the cones $\sigma_v$, that
is, they glue together into a continuous function on $\ct$. We
conclude that $H^*_T(X)$ is identified with the algebra of continuous
functions on $\ct$ which are piecewise polynomial with respect to
the subdvision $(\sigma_F)$ ($F$ a face of $P$). Moreover, the image
of $S$ in $H^*_T(X)$ is identified to the algebra of polynomial
functions. We refer to \cite{Br2} for more on the relations between
continuous, piecewise polynomial functions and toric manifolds.
\bigskip

We now generalize this description to the class of all compact
multiplicity-free Hamiltonian spaces, in the following sense.
\bigskip

\noindent
{\bf Definition} \cite{GS}, \cite{W} A Hamiltonian $G$-space $X$ is
{\sl multiplicity-free} if $X$ is connected and the preimage under the
moment map of any coadjoint $G$-orbit consists of finitely many
$G$-orbits.
\bigskip

If moreover $X$ is compact, then the fibers of the moment map are
connected (see \cite{Lerman} for a simple proof of this fact). Thus,
multiplicity-free amounts to: the preimage of each orbit under the
moment map  is an unique orbit.

By \cite{Br1}, a complex projective $G$-manifold $X$ is
multiplicity-free if and only if it is {\sl spherical} for the action
of the  complexification $G^{\C}$, that is, a Borel subgroup of
$G^{\C}$ has a dense orbit in $X$.

For compact multiplicity-free spaces, we have the following sharper
version of Lemma \ref{fix}. 
\begin{lemma}\label{fixbis}
Let $X$ be a compact multiplicity-free $G$-space and $\Gamma\subset G$
a closed subgroup. Then $X^{\Gamma}$ has only finitely many components,
and each of them is a compact multiplicity-free $G^{\Gamma}$-space. In
particular, the fixed point set of a maximal torus of $G$ is finite.
Moreover, denoting by $X_1$ the set of all $x\in X$ such that the rank
of the isotropy group $G_x$ is at least $rk(G)-1$, and identifying
$(\g^*)^T$ with $\ct^*$, the intersection $\mu(X_1)\cap\ct^*$ is a
finite union of segments with ends in $\mu(X^T)$.
\end{lemma}
\begin{proof}
Let $Y$ be a component of $X^{\Gamma}$ and let $x\in Y$. Then we know
that $\mu^{-1}(\mu(Gx))=Gx$. Thus, $\mu\vert_Y^{-1}(\mu(G^{\Gamma}x))$
is contained in $(Gx)^{\Gamma}$. But the set $(Gx)^{\Gamma}$ is a
finite union of orbits of $G^{\Gamma}$. It follows that $Y$ is
multiplicity-free.

In particular, each component $Y$ of $X^T$ is a multiplicity-free
space for the action of $G^T=T$. Thus, $Y$ must be a point.

For the last assertion, observe first that $X_1$ is the union of the
sets $GX^{T'}$ where $T'\subset T$ is a subtorus of codimension 1.
Moreover, the number of subsets $X^{T'}$ is finite (this can be seen
by linearizing the action of $T$ around fixed points). Thus, it is
enough to check that each $\mu(GX^{T'})\cap\ct^*$ is a finite union of
segments with ends in $\mu(X^T)$. But
$$
\mu(GX^{T'})\cap\ct^*=G\mu(X^{T'})\cap\ct^*
$$
and $\mu(X^{T'})$ is contained in $(\g^*)^{T'}=G^{T'}\ct^*$, whence
$$
\mu(GX^{T'})\cap\ct^*=W(\mu(X^{T'})\cap\ct^*).
$$
So we can replace $(G,X)$ by $(G^{T'}/T',Y)$ where $Y$ is a component
of $X^{T'}$, and thus we can assume that the rank of $G$ is one. Then
$\mu(X)\cap\ct^*$ is a segment, or a union of two segments with ends
in $\mu(X^T)$ \cite{I}.
\end{proof}
\bigskip

Observe that $W$ acts on $X^T$ and that both $\mu(X^T)$ and
$\mu(X_1)\cap \ct^*$ are $W$-invariant subsets of $\ct^*$. Choose a
Weyl chamber $\ct^*_+\subset \ct^*$ and set
$$
\mu(X^T)\cap \ct^*_+=\{\lambda_1,\ldots,\lambda_m\}.
$$
Then $\mu(X^T)=W\lambda_1\cup\cdots\cup W\lambda_m$. For 
$1\leq j\leq m$, choose $x_j\in X^T$ such that $\mu(x_j)=\lambda_j$
and denote by $W_j$ the isotropy group of $x_j$ in $W$. Then the 
$W$-orbits in $X^T$ are the cosets $W/W_j$ ($1\leq j\leq m$). As
a consequence, we have an isomorphism (by Proposition \ref{deg})
$$
H^*_T(X^T)^W\simeq\prod_{j=1}^m S^{W_j}.
$$
We now describe the image of $H^*_G(X)\simeq H^*_T(X)^W$ in 
$H^*_T(X^T)^W$ under restriction to fixed points $i_T^*$.
\begin{theorem}\label{multfree}
Let $X$ be a compact multiplicity-free space under a compact 
connected Lie group $G$. Then, with notation as above, the algebra
$H^*_G(X)$ is isomorphic via $i_T^*$ to the subalgebra of $S^m$
consisting of all $m$-uples $(f_1,\ldots,f_m)$
such that: 

1) each $f_j$ is in $S^{W_j}$, and

2) $f_j\equiv w(f_k)$ (mod $\lambda_j-w(\lambda_k)$) whenever $w\in W$
and the segment $[\lambda_j,w(\lambda_k)]$ is a component of 
$\mu(X_1)\cap\ct^*$.

Moreover, $i_T^*(S^W)$ consists of all tuples $(f,\ldots,f)$ where
$f\in S^W$.
\end{theorem}
\begin{proof}
Let $\Gamma\subset T$ be a subtorus of codimension one. Let $\chi$ be
a character of $T$ with kernel $\Gamma$ and let $Y$ be a connected
component of $X^{\Gamma}$. By Lemma \ref{fixbis}, the
$G^{\Gamma}$-variety $Y$ is multiplicity-free. Set
$G^{\Gamma}/\Gamma:=H$; then $H$ is a compact connected Lie group of
rank 1, and thus it is isomorphic to $S^1$, ${\rm SU}(2)$
or ${\rm SO}(3)$. Therefore, the $H$-multiplicity free space $Y$ must
be of dimension at most 4. Two cases can occur:

\noindent
1) $Y$ is two-dimensional. Then $Y$ is isomorphic to complex projective
line, and $Y^T$ consists of two fixed points $y$, $z$.
Restriction to these fixed points identifies $H^*_T(Y)$ to
the set of all $(f_y,f_z)\in S\times S$ such that
$$
f_y\equiv f_z~({\rm mod}\,\chi)~.\leqno(1)
$$ 

If $\chi$ is not a root of $(G,T)$, then $H\simeq S^1$ and $\mu(Y)$
is the segment $[\mu(y),\mu(z)]$ in $\ct^*$. Thus, this segment lies in
$\mu(X_1)\cap\ct^*$. On the other hand, if $\chi$ is a root,
let $s\in W$ be the corresponding reflection and let 
$\tilde s$ be a representative of $s$ in the normalizer of $T$. Then
$Y$ is invariant under $\tilde s$, and  
$$
H^*_T(Y)^s=\{(f_y,s(f_y))~\vert~f_y\in S\}
$$
(indeed, $f-s(f)$ is divisible by $\chi$ for any $f\in S$).
\noindent

2) $Y$ is four-dimensional. Then $H$ is isomorphic to ${\rm SU}(2)$
or ${\rm SO}(3)$, and the $H$-variety $Y$ is either a
rational ruled surface, or the projectivization of a three-dimensional
complex representation of ${\rm SU}(2)$ (see \cite{I} and \cite{Au}
Chapter IV, Appendix A). In the former case, $Y^T$ consists of four
points $y$, $s(y)$, $z$, $s(z)$ where $s$ is the non-trivial element
of the Weyl group of $(H,T/\Gamma)$; we may assume that the segment 
$[\mu(y),\mu(z)]$ lies in $\mu(Y)\cap\ct^*$. It is
easy to check that restriction to fixed points maps $H^*_T(Y)$ onto
the set of all quadruples $(f_y,f_{s(y)},f_z,f_{s(z)})\in S^4$ such
that $$
f_y\equiv f_{s(y)}\equiv f_z\equiv f_{s(z)}~({\rm mod}\,\chi),~
f_y+f_{s(y)}\equiv f_z+f_{s(z)}~({\rm mod}\,\chi^2)~.\leqno(2)
$$
It follows that
$$
H^*_T(Y)^s=\{(f_y,f_z)\in S\times S~\vert~
f_y\equiv f_z~ ({\rm mod}\,\chi)\}.
$$

In the latter case, we have similarly $Y^T=\{y,s(y),z\}$ where
$z=s(z)$, and $H^*_T(Y)$ consists of the triples
$(f_y,f_{s(y)}, f_z)\in S^3$ such that
$$f_y\equiv f_{s(y)}\equiv f_z~({\rm mod}\,\chi),~
f_y+f_{s(y)}\equiv 2f_z~({\rm mod}\,\chi^2)~.\leqno(3)$$
It follows that
$$
H^*_T(Y)^s=\{(f_y,f_z)\in S\times S^s~\vert~
f_y\equiv f_z~({\rm mod}\,\chi)\}~.
$$
We conclude that the image of $i_T^*$ is defined by our congruences.
Observe that the image of 
$i_T^*:H^*_T(X)\to H^*_T(X^T)$ is defined by congruences of the
form (1), (2) or (3).
\end{proof}
\bigskip

\noindent
{\bf Examples} 

1) (coadjoint orbits) Let $X$ be the $G$-orbit of
$\lambda\in\g^*$; we may assume that $\lambda\in\ct^*_+$. Then
$\mu:X\to\g^*$ is the inclusion map, whence 
$\mu(X^T)\cap\ct^*=\{\lambda\}$ and 
$\mu(X_1)\cap\ct^*=W\cdot\lambda$. Then Theorem \ref{multfree}
reduces to the isomorphism
$$
H^*_G(G\lambda)=S^{W_{\lambda}}~
$$
which follows more directly from the isomorphism
$$
H^*_G(G\lambda)=H^*_G(G/G_{\lambda})=H^*_{G_{\lambda}}(pt).
$$

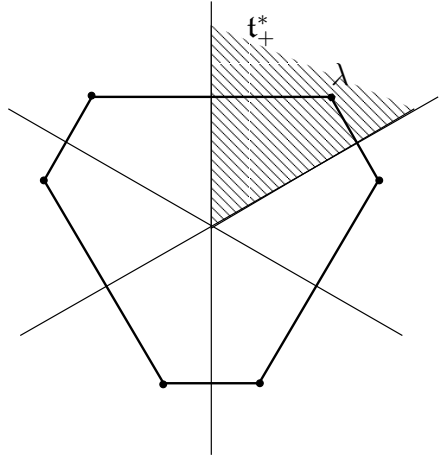
\begin{center}

\begin{picture}(0,0)%
\epsfig{file=figure2.pstex}%
\end{picture}%
\setlength{\unitlength}{0.00083300in}%
\begingroup\makeatletter\ifx\SetFigFont\undefined%
\gdef\SetFigFont#1#2#3#4#5{%
  \reset@font\fontsize{#1}{#2pt}%
  \fontfamily{#3}\fontseries{#4}\fontshape{#5}%
  \selectfont}%
\fi\endgroup%
\begin{picture}(2724,2874)(814,-2098)
\put(2851,239){\makebox(0,0)[lb]{\smash{\SetFigFont{12}{14.4}{\rmdefault}{\mddefault}{\updefault}$\lambda$}}}
\put(2326,539){\makebox(0,0)[lb]{\smash{\SetFigFont{12}{14.4}{\rmdefault}{\mddefault}{\updefault}$\mathfrak t^*_+$}}}
\end{picture}

\end{center}

\medskip

2) (toric varieties). If $T$ is a torus and $X$ a projective toric
manifold for the complexification $T^{\C}$, then the set
$\mu(X_1)\cap\ct^*$ is the union of all edges of the polytope
$\mu(X)$. So Theorem \ref{multfree} gives back the description of
$H^*_T(X)$ found in the second example of this section.

3) (complete conics). Let $V$ be the vector space of
quadratic forms on $\C^3$, let $V^*$ be the dual space,
and let $\CP=\CP(V)\times\CP(V^*)$ be the product of their
projectivizations. Let $X\subset\CP$ be the closure of the set of
classes $([A],[A^{-1}])$ where $A\in V$ is non-degenerate and
$A^{-1}\in V^*$ is the dual quadratic form. Then
$X$ is a complex projective manifold, called the {\sl space of
complete conics}. Moreover, $X$ is spherical for the natural action of
${\rm GL}(3)$, and hence multiplicity-free for the action of the
maximal compact subgroup $G:=U(3)$. 

\newpage

The set $\mu(X_1)\cap\ct^*$ is given by the following picture.

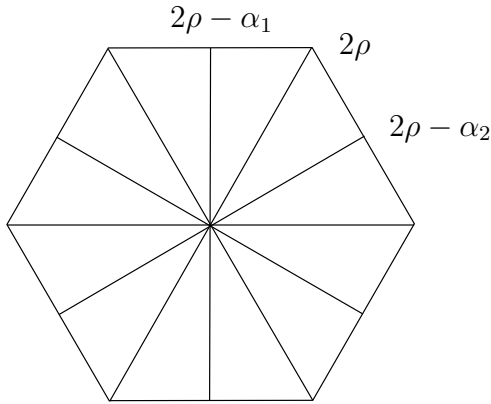
\begin{center}
\begin{picture}(0,0)%
\epsfig{file=figure3.pstex}%
\end{picture}%
\setlength{\unitlength}{0.00083300in}%
\begingroup\makeatletter\ifx\SetFigFont\undefined%
\gdef\SetFigFont#1#2#3#4#5{%
  \reset@font\fontsize{#1}{#2pt}%
  \fontfamily{#3}\fontseries{#4}\fontshape{#5}%
  \selectfont}%
\fi\endgroup%
\begin{picture}(2583,2505)(1265,-1975)
\put(2311,374){\makebox(0,0)[lb]{\smash{\SetFigFont{12}{14.4}{\rmdefault}{\mddefault}{\updefault}$2\rho-\alpha_1$}}}
\put(3361,209){\makebox(0,0)[lb]{\smash{\SetFigFont{12}{14.4}{\rmdefault}{\mddefault}{\updefault}$2\rho$}}}
\put(3676,-301){\makebox(0,0)[lb]{\smash{\SetFigFont{12}{14.4}{\rmdefault}{\mddefault}{\updefault}$2\rho-\alpha_2$}}}
\end{picture}
\end{center}

It follows that the algebra $H^*_G(X)$ consists of all triples
$(f,f_1,f_2)$ in $S\times S^{s_1}\times S^{s_2}$ such that
$f\equiv f_1$ (mod $\alpha_1$), $f\equiv f_2$ (mod $\alpha_2$) and
that $f_1\equiv s_{\alpha_1+\alpha_2}(f_2)$ (mod $2\alpha_1+\alpha_2$)
where $\alpha_1$, $\alpha_2$ are the simple roots, with corresponding
reflections $s_1$, $s_2$.

The variety of complete conics is the basic example of a ``complete
symmetric variety'', a nice compactification of a symmetric space
\cite{DP}. The cohomology ring of complete symmetric varieties is
described by other methods in \cite{BDP} and \cite{LP}.

\section{Equivariant Chow groups}

We begin by recalling the definition and some basic features of
(usual) Chow groups; a complete exposition can be found in
\cite{Fulton}.

Let $X$ be a scheme of finite type over the field of complex numbers
(most of what follows holds more generally over any algebraically
closed field). The {\sl group of algebraic cycles} on $X$ is the
abelian group $Z_*(X)$ freely generated by symbols $[Y]$ where
$Y\subset X$ is a subvariety (that is, a closed subscheme which is a
variety). Observe that the group $Z_*(X)$ is graded by dimension. 

For $Y$ as before, and for a rational function $f$ on $Y$, we define
the divisor of $f$, an element of $Z_*(X)$, by
$$
div(f):=\sum_D ord_D(f)[D]
$$

\noindent (sum over all prime divisors $D\subset Y$), where $ord_D(f)$ denotes
the order of the
\newpage

\noindent zero or pole of $f$ along $D$. The cycles $div(f)$
for $Y$ and $f$ as before, generate a graded subgroup 
$Rat_*(X)\subset Z_*(X)$: the {\sl group of rationally trivial
cycles}.

By definition, the {\sl Chow group} $A_*(X)$ is the quotient
$Z_*(X)/Rat_*(X)$; it is a graded abelian group with homogeneous
components $A_n(X)$ where $0\leq n\leq \dim(X)$. The top degree
component $A_{\dim(X)}(X)$ is freely generated by the top dimensional
irreducible components of $X$.

If $X$ is equidimensional (that is, if all irreducible
components of $X$ have the same dimension), then we set
$A^n(X):=A_{\dim(X)-n}(X)$, the group of cycles of codimension $n$,
and we denote by $A^*(X)$ the Chow group graded by codimension. If
moreover $X$ is smooth, then there is an intersection product on
$A^*(X)$ which makes it a graded ring.

We now review the functorial properties of Chow groups.
Any proper morphism $\pi:X\to Y$ defines a push-forward homomorphism
$\pi_*:A_*(X)\to A_*(Y)$ which preserves degree. On the
other hand, a flat morphism $\pi:X\to Y$ with $d$-dimensional fibers
induces a pull-back homomorphism $\pi^*:A_*(Y)\to A_*(X)$ which
increases degree by $d$. If moreover all fibers are isomorphic to
affine space, then $\pi^*$ is an isomorphism \cite{Gil}. In
particular, the projection of a vector bundle induces isomorphisms on
Chow groups. This can be seen as a version of homotopy invariance.
Finally, any local complete intersection morphism $\pi:X\to Y$ of
codimension $d$ (in the sense of \cite{Fulton} 6.6) induces a
homomorphism $\pi^*:A_*(Y)\to A_*(X)$ which decreases degree by $d$.
In particular, if $X$ and $Y$ are smooth, then any morphism
$\pi:X\to Y$ defines a ring homomorphism $\pi^*:A^*(Y)\to A^*(X)$
which preserves degree.

For a scheme $X$ and a closed subscheme $Y$, we have a short exact
sequence
$$
A_*(Y)\to A_*(X)\to A_*(X\setminus Y)\to 0.
$$
It can be extended to a long exact sequence by introducing higher
Chow groups (which will not be considered here).

Any vector bundle $E$ on $X$ has Chern classes $c_j(E)$, which are
homogeneous operators on $A_*(X)$ of respective degree $-j$. If $X$ is
smooth and $i:Y\to X$ denotes the inclusion of a smooth subvariety of
codimension $d$, then the composition $i^*\circ i_*$ is multiplication
by $c_d(N_{X,Y})$, the top Chern class of the normal bundle to $Y$ in
$X$ (this holds more generally if $Y$ is a local complete intersection
in $X$).

Each subvariety $Y\subset X$ has a homology class $cl_X(Y)$ of degree
$2\dim(Y)$ in Borel-Moore homology $H_*(X,\Z)$ \cite{Fulton} 19.1. The
assignement $Y\mapsto cl_X(Y)$ defines a {\sl cycle map} 
$$
cl_X:A_*(X)\to H_*(X,\Z)
$$
which doubles degrees. If moreover $X$ is smooth of (complex)
dimension $N$, then $H_j(X,\Z)=H^{2N-j}(X,\Z)$ and thus, we have
a cycle map
$$
cl_X:A^*(X)\to H^*(X,\Z),
$$
which is a (degree doubling) ring homomorphism. The cycle map is an
isomorphism if $X$ has a cellular decomposition \cite{Fulton} Example
19.1.11. But in general, the cycle map is far from being injective or
surjective.
\bigskip

\noindent
{\bf Examples} 

1) (curves) Let $X$ be a smooth projective curve of
genus $g$. Then $A_1(X)\simeq\Z$ is freely generated by $[X]$, and is
isomorphic to $H_2(X,\Z)$. The cycle map $A_0(X)\to H_0(X,\Z)$ is the
degree map; its kernel is the Jacobian variety of $X$, a compact torus
of dimension $2g$. The cokernel of the cycle map $A_*(X)\to H_*(X,\Z)$
is the group $H_1(X,\Z)$, isomorphic to $\Z^{2g}$.

2) (linear algebraic groups) Let $G$ be a connected linear algebraic
group and $B\subset G$ be a Borel subgroup. Then the flag manifold
$G/B$ has a cellular decomposition by Schubert cells, and thus its
Chow ring is isomorphic to its integral cohomology ring. Let 
$T\subset B$ be a maximal torus and let $\Xi(T)$ be its character
group. Then we have the characteristic homomorphism $c:\Xi(T)\to
A^1(G/B)$ which maps each character to the Chern class of the
associated line bundle over $G/B$. It extends to an algebra
homomorphism
$$
c:S_{\Z}\to A^*(G/B)
$$
where $S_{\Z}$ is the symmetric algebra of $\Xi(T)$ over the integers.

We can choose a maximal compact subgroup $G_c\subset G$ such that the
group $T\cap G_c=:T_c$ is a maximal compact torus of $T$; then the map
$G_c/T_c\to G/B$ is an homeomorphism. Now, by Proposition \ref{deg},
the characteristic homomorphism induces an isomorphism
$$
S/(S^W_+)\to A^*(G/B)_{\Q}.
$$

Moreover, the Chow ring of $G$ is isomorphic to the quotient of
$A^*(G/B)$ by its ideal generated by $c(\Xi(T))$ \cite{Gro}
p. 21. To see this, choose a basis $(\chi_1,\ldots,\chi_{\ell})$ of
$\Xi(T)$. Consider the action of $T$ on $\C^{\ell}$ with
weights $\chi_1,\ldots,\chi_{\ell}$. Then $T$ embeds into $\C^{\ell}$
as the complement of the coordinate hyperplanes. Let
$E:=G\times_T\C^{\ell}$ be the associated vector bundle over $G/T$.
Then $E$ is a direct sum of line bundles $L_1,\ldots,L_{\ell}$, and
$G=G\times_T T$ embeds into $E$ as the complement of the union of zero
sections $E^{(j)}:=\oplus_{i\neq j}L_i$ ($1\leq j\leq \ell$).
Using the exact sequence
$$
\oplus_{j=1}^{\ell} A^*(E^{(j)})\to A^*(E)\to A^*(G)
$$
and the fact that the image of $A^*(E^{(j)})$ in $A^*(E)$ is the ideal
generated by $[E^{(j)}]$, we see that
$$\displaylines{
A^*(G)=A^*(E)/([E^{(1)}],\ldots,[E^{(\ell)}])
=A^*(G/T)/(c(\chi_1),\ldots,c(\chi_{\ell}))
\hfill\cr\hfill
=A^*(G/B)/(c(\Xi(T))).
\cr}$$
In particular, because $c$ is surjective over $\Q$, the positive
degree part of $A^*(G)$ is finite.

On the other hand, the rational cohomology algebra $H^*(G)$ is a free
exterior algebra on $\ell$ generators \cite{Hopf}. Thus, the cycle map
$cl_G$ is not an isomorphism over $\Q$, except when $G$ is unipotent.
\bigskip

Now we introduce the equivariant Chow groups, after Edidin and Graham
\cite{EG2}. Let $G$ be a linear algebraic group and let $n$ be a non
negative integer. As in Example 3 of the first section, we can find a 
$G$-module $V$ and a $G$-invariant open subset $U\subset V$ (depending
on $n$) satisfying the following conditions.

1) The quotient $U\to U/G$ exists and is a principal $G$-bundle.

2) The codimension of $V\setminus U$ in $V$ is larger than $n$.

Observe that $U\to U/G$ is an approximation of the universal
$G$-bundle by an algebraic bundle.

Now let $X$ be a scheme with a $G$-action such that $X$ can be covered
by invariant quasi-projective open subsets (by \cite{S}, this
assumption is fulfilled when $X$ is normal); we will say that $X$ is a
$G$-scheme. Then the quotient of $X\times U$ by the diagonal action of
$G$ exists as a scheme; we denote this quotient by $X\times_G U$. 

If $X$ is equidimensional, we define its {\sl equivariant Chow group
of degree $n$} by
$$
A^n_G(X):=A^n(X\times_G U).
$$
This makes sense because $X\times_G U$ is equidimensional, and also
because this does not depend on the choice of $U$. Indeed, if
$U'\subset V'$ is another choice, then the quotient of 
$X\times U\times V'$ by the diagonal action of $G$ exists and defines
a map
$$
p:(X\times U\times V')/G\to X\times_G U.
$$
Observe that $p$ is smooth with fibers isomorphic to $V'$. Thus, we
have isomorphisms 
$$
A^n(X\times_G U)\simeq A^n((X\times U\times V')/G)\simeq 
A^n((X\times U\times U')/G),
$$
the latter being a consequence of assumption (2) for $U'$. It follows
that $A^n(X\times_G U)\simeq A^n(X\times_G U')$.

For $X$ not necessarily equidimensional, we set
$$
A_n^G(X):=A_{n-\dim(G)+\dim(U)}(X\times_G U).
$$
Arguing as above, we see that this group is independent of $U$, so
that we can define the {\sl equivariant Chow group}
$$
A_*^G(X)=\bigoplus_{n\in\Z} \,A_n^G(X).
$$
By definition, we have $A_n^G(X)=0$ for $n>\dim(X)$; but $A_n^G(X)$
may be non trivial for $n<0$. For equidimensional $X$, we have
$$
A^n_G(X)=A_{\dim(X)-n}^G(X).
$$

We now list some properties of equivariant Chow groups; let us mention
first that they satisfy the functorial properties of usual Chow
groups, and that each linearized vector bundle has equivariant Chern
classes.

If $X$ is smooth, then the same holds for each $X\times_G U$. It
follows that $A^*_G(X)$ is a graded ring for the intersection product.
In particular, $A^*_G(pt)$ is a graded ring. 

For arbitrary $X$, the projection $X\times U\to U$ descends to a flat
map
$$
p_X:X\times_G U\to U/G.
$$
Therefore, $A_*^G(X)$ is a graded $A^*_G(pt)$-module.
The image of $p_X$ is the smooth variety $U/G$, and the fibers are
isomorphic to $X$. Moreover, any two points in $U/G$  can be joined by
a chain of rational curves (because the same holds in $U$, an open
subset of a linear space). Thus, pull-back to a fiber is a
well-defined map $A_*^G(X)\to A_*(X)$ invariant under the action of
$A^*_G(pt)$. So we obtain a map
$$
\rho:A_*^G(X)/(A^+_G(pt))\to A_*(X).
$$

Any $G$-invariant subvariety $Y\subset X$ defines an equivariant class
$[Y]$ in $A_*^G(X)$: indeed, set 
$$
[Y]:=[Y\times_G U]\in A_*(X\times_G U)
$$
for $U$ as above such that
${\rm codim}_V(V\setminus U)>{\rm codim}_X(Y)$. The image of $[Y]$ in
$A_*(X)$ is the (usual) class of $Y$. In fact, we have a degree
doubling cycle map
$$
cl_X^G:A_*^G(X)\to H_*^G(X,\Z)
$$
to equivariant Borel-Moore homology \cite{EG2}. For $X$ smooth,
we have a cycle map 
$$
cl_X^G:A^*_G(X)\to H^*_G(X,\Z)
$$
to equivariant cohomology, which lifts the usual cycle map $cl_X$.
 
If $G$ acts on $X$ with a quotient $X\to X/G$ which is a principal
$G$-bundle, then we have a smooth map $X\times_G V\to X/G$ with fibers
isomorphic to $V$. It follows that 
$$
A_{n-\dim(G)+\dim(U)}(X\times_G U)=A_{n-\dim(G)+\dim(V)}(X\times_G V)
=A_{n-\dim(G)}(X/G).
$$
Thus, $A_*^G(X)$ is isomorphic to $A_*(X/G)$ (with degree shifted by 
$\dim(G)$). A much deeper result is due to Edidin and Graham: they
proved that $A_*^G(X)$ is isomorphic to $A_*(X/G)$ over $\Q$,
whenever $G$ acts on $X$ with finite isotropy groups \cite{EG2}. 

Let $X$ be a $G$-scheme and let $H\subset G$ be a closed subgroup.
Then, for $U$ as above, the quotient $U\to U/H$ exists and is a
principal $H$-bundle. Thus, there is a smooth map 
$X\times_H U\to X\times_G U$ with fiber $G/H$, which induces a map
$$
A_*^G(X)\to A_*^H(X)
$$
of degree $\dim(G/H)$. If moreover $H$ is a Levi subgroup of $G$, then
$G/H$ is the unipotent radical of $G$ and hence is isomorphic to
affine space. Thus, $A_*^G(X)$ is isomorphic to $A_*^H(X)$.

Our latter remark reduces many questions on equivariant Chow groups to
the case of reductive groups. From now on, we assume that $G$ is
reductive and connected; we denote by $T\subset G$ a maximal torus, by
$W$ its Weyl group and by $B$ a Borel subgroup of $G$ containing $T$.
Let $S_{\Z}$ (resp. $S$) be the symmetric algebra
over the integers (resp. over the rationals) of the character group
$\Xi(T)\simeq\Xi(B)$. Then we have the following analogue of
Proposition \ref{deg}, due to Edidin and Graham \cite{EG1}, \cite{EG2}.

\begin{theorem}\label{isom}
Notation being as above, the graded ring $A^*_T(pt)$ is isomorphic to
$S_{\Z}$. Moreover, the map $A^*_G(pt)\to A^*_T(pt)$ is injective over
$\Q$ and identifies $A^*_G(pt)_{\Q}$ to $S^W$. Finally, for any
$G$-scheme $X$, we have isomorphisms 
$A_*^G(X)_{\Q}\simeq A_*^T(X)^W_{\Q}$ and
$S\otimes_{S^W}A_*^G(X)_{\Q}\simeq A_*^T(X)_{\Q}$.
\end{theorem}

Indeed, for any $\chi\in\Xi(T)$, consider the line bundle 
$L(\chi)=U\times_T \C\chi$ on $U/T$ as in Example 2 of Section 1. The
first Chern class of $L(\chi)$ defines an element $c(\chi)$ of
$A^1_T(pt)$, and the assignement $\chi\mapsto c(\chi)$ extends to the
characteristic homomorphism $c:S_{\Z}\to A^*_T(pt)$. Using an explicit
description of $U\to U/T$ as in Example 1 of Section 1, we see that
$c$ is an isomorphism.

To prove the other statements, one begins by observing that the map 
$$
X\times_T U\to X\times_G U
$$ 
factors through $X\times_B U$. Moreover, the map 
$$
X\times_T U\to X\times_B U
$$ 
is smooth with fiber $B/T$ isomorphic to affine space, so it induces
an isomorphism of Chow groups. Finally, the map 
$$
X\times_B U\to X\times_G U
$$
is smooth and proper with fiber the flag variety $G/B$. Then an
analogue of the Leray-Hirsch theorem holds for this map \cite{EG1},
and we conclude as in the proof of Proposition \ref{deg}.
\bigskip

For a $T$-scheme $X$, the group $A_*^T(X)$ is a $S_{\Z}$-module
containing the equivariant classes $[Y]$ where $Y\subset X$ is a
$T$-stable subvariety. If moreover $f$ is a rational function on $Y$
which is an eigenvector of $T$ of weight $\chi$, then the support of
its divisor is $T$-invariant, and thus we can define a class $div(f)$
in $A_*^T(X)$. In fact, we have 
$$
div(f)=\chi[Y].
$$
Indeed, $f$ can be seen as a rational section
of the pull-back of the line bundle $L(\chi)=U\times_T \C\chi$ to 
$U\times_T Y$. Thus, $div(f)$ represents the pull-back of $c(\chi)$
evaluated on the class of $U\times_T Y$, that is, $\chi[Y]$.

It turns out that the $S_{\Z}$-module $A_*^T(X)$ is defined by
generators $[Y]$ and relations $[div(f)]-\chi[Y]$ as above \cite{Br3}
2.1. The proof is based on an explicit construction of $U$ and $U/T$
as toric varieties, combined with a result of Fulton, MacPherson,
Sottile and Sturmfels: For any scheme $X$ with an action of a
connected solvable linear algebraic group $\Gamma$, the usual Chow
group $A_*(X)$ is defined by generators $[Y]$ (where $Y\subset X$ is a
$\Gamma$-stable subvariety) and relations $[div(f)]$ (where $f$ is a
rational function, eigenvector of $\Gamma$, on $Y$ as before)
\cite{FMSS}. As a consequence, we obtain the following analogue of
Proposition \ref{kir}, valid for equivariant Chow groups of 
{\sl arbitrary} schemes \cite{Br3} 2.3.

\begin{proposition}\label{iso}
Let $G$ be a connected reductive algebraic group, $T\subset G$ a
maximal torus, and $X$ a $G$-scheme. Then the maps
$$
\rho_T:A_*^T(X)/(S_{\Z,+})\to A_*(X)
$$
and
$$
\rho_G: A_*^G(X)_{\Q}/(S^W_+)\to A_*(X)_{\Q}
$$
are isomorphisms.
\end{proposition}
\begin{proof} The first isomorphism follows immediatly from the
results above; it can also be proved directly, as follows. Fixing a
degree $n$ and using the isomorphism 
$A_n(X)\simeq A_{n+\dim(U)}(X\times U)$, we reduce to checking that
the pull-back map 
$$
A_*((X\times U)/T)/(S_{\Z,+})\to A_*(X\times U)
$$ 
is an isomorphism. But this follows from the argument of Example 2
above, where it was shown that $A_*(G/T)/(S_{\Z,+})\to A_*(G)$ is an
isomorphism. 

The second isomorphism follows from the first one, combined with the
isomorphism $S\otimes_{S^W}A_*^G(X)_{\Q}\to A_*^T(X)_{\Q}$ of Theorem
\ref{isom}.
\end{proof}
\bigskip

Consider for example $X=G$ where $G$ acts by left multiplication. Then
$A^*_T(G)=A^*(G/T)=A^*(G/B)$ and we recover the structure of $A^*(G)$
given in Example 2 above. For $X=G/B$, we have $A^*_G(X)_{\Q}=S$ and
we recover the description of $A^*(G/B)_{\Q}$ obtained earlier by
comparing with cohomology. 

More generally, let us apply this theory to a description of rational
Chow rings of homogeneous spaces $G/H$ where $G$ is a connected linear
algebraic group and $H\subset G$ is a closed connected subgroup. As
above, we may assume that $G$ and $H$ are reductive. 

Let $T_H$ be a maximal torus of $H$ with Weyl group $W_H$, and let
$S_H$ be the symmetric algebra over $\Q$ of the character group of
$T_H$. Finally, let $T$ be a maximal torus of $G$ containing $T_H$ and
let $W$ be the Weyl group of $T$; then we have a restriction map 
$S\to S_H$. Observe that this maps sends $S^W$ to $S_H^{W_H}$. Indeed,
the complexified space $S_{\C}$ can be identified to the algebra of
polynomial functions on $\ct$, and the subspace $S^W_{\C}$ is then
spanned by restrictions to $\ct$ of characters of $\g$-modules
\cite{Bour} VIII.3.3. The latter restrict to characters of
$\h$-modules, which are $W_H$-invariants. We denote by $(S^W_+)$ the
ideal of $S_H^{W_H}$ generated by restrictions of homogeneous elements
of positive degree of $S^W$. 
\begin{corollary}\label{hom} 1) Notation being as above, we have
isomorphisms of graded rings  
$$
A^*(G/H)_{\Q}\simeq S_H^{W_H}/(S^W_+)
$$
and
$$
A^*_T(G/H)_{\Q}\simeq S\otimes_{S^W} S_H^{W_H}.
$$

\noindent
2) The following conditions are equivalent:

(i) The rational cohomology of $G/H$ vanishes in odd degree.

(ii) The ranks of $G$ and $H$ are equal.
 
(iii) The cycle map $cl_{G/H}:A^*(G/H)_{\Q}\to H^*(G/H)$ is an
isomorphism. 
\end{corollary}
\begin{proof} 1) For the first isomorphism, observe that
$A^*_G(G/H)=A^*_H(pt)$ and apply Proposition \ref{iso}. For the second
one, recall that $A^*_T(G/H)_{\Q}$ is isomorphic to 
$S\otimes_{S^W} A^*_G(G/H)_{\Q}$.

2) (i)$\Rightarrow$(ii) If $rk(H)<rk(G)$ then $T$ acts on $G/H$
without fixed point. It follows that the topological Euler
characteristic of $G/H$ is zero. Thus, $G/H$ must have non trivial odd
cohomology.

(ii)$\Rightarrow$(iii) By assumption, we have $T_H=T$ and therefore,
$$
H^*(G/H)=H^*_H(G)=H^*_T(G)^{W_H}=H^*(G/T)^{W_H}
$$
with similar isomorphisms for (equivariant) Chow rings. Now the cycle map
$$
A^*(G/T)=A^*(G/B)\to H^*(G/B)=H^*(G/T)
$$
is a $W$-equivariant isomorphism. This implies our statement.

(iii)$\Rightarrow$(i) is obvious.
\end{proof}
\bigskip

Recall that any homogeneous space $G/H$ under a connected linear
algebraic group has the structure of a vector bundle over the
homogeneous space $G_c/H_c$ where $H_c$ is a maximal compact subgroup
of $H$, and $G_c\supset H_c$ is a maximal compact subgroup of $G$. In
particular, $G/H$ has the same cohomology as $G_c/H_c$. The cohomology
ring of the latter has been intensively studied, see e.g.
\cite{Borel}; its structure is in general more complicated than that
of $A^*(G/H)_{\Q}$. Concerning the {\sl integral} Chow ring
$A^*(G/H)$, we are led to the following
\bigskip

\noindent
{\bf Question 1.} How to describe torsion in Chow rings of
homogeneous spaces? 
\bigskip

As seen in Example 2 above, the Chow ring of a connnected linear
algebraic group is torsion-free (that is, isomorphic to $\Z$) if and
only if the characteristic homomorphism $c:S_{\Z}\to A^*(G/B)$ is
surjective, that is, the group is special in the sense of \cite{Gro}. 
\bigskip

\noindent
{\bf Question 2.} Notation being as above, when is $A^*(G/H)$
finite in positive degree? Equivalently, when is the restriction map 
$S^W\to S_H^{W_H}$ surjective?  Or, denoting by $\C[\g]$ the ring of
polynomial functions on $\g$ and using the Chevalley restriction
theorem \cite{Bour} VIII.3.3: when is the restriction map 
$\C[\g]^G\to \C[\h]^H$ 
surjective? 
\bigskip

This holds e.g. when $H$ is isomorphic to ${\rm SL}(2)$ or to 
${\rm PSL}(2)$: indeed, the Killing form of $\g$ restricts then to a
non zero quadratic element of $\C[\h]^H$, which generates this
algebra. But this does not hold if $H$ is a proper subgroup
of maximal rank of $G$, because we then have $T_H=T$ and $W_H\neq W$.
It turns out that this does not hold for certain pairs of semisimple
groups $(G,H)$ as well. The following example was pointed out by Bram
Broer: let $G$ be semisimple of type $E_6$ and $H\subset G$ be a
maximal semisimple subgroup of a parabolic subgroup of type $D_5$.
Then the sequence of degrees of a minimal system of homogeneous
generators for $S^W$ (resp. $S_H^{W_H}$) is $(2,5,6,8,9,12)$ (resp.
$(2,4,5,6,8)$): a generator of degree 4 of $S_H^{W_H}$ cannot be in
the image of $S^W$.

\section{Localization and equivariant multiplicities}

As in the previous section, we consider schemes over the field of
complex numbers. If $T$ is a torus and $X$ a $T$-scheme, we denote by
$$
i_T:X^T\to X
$$ the inclusion of the fixed point subscheme. Then $i_T$ induces a
homomorphism of equivariant Chow groups
$$
i_{T*}:A_*^T(X^T)\to A_*^T(X^T)
$$
which is $S_{\Z}$-linear (here $S_{\Z}$ denotes as before the
symmetric algebra over the integers of the character group of $T$).
Moreover, we have a natural isomorphism 
$$
A_*^T(X^T)\simeq S_{\Z}\otimes A_*(X^T).
$$
Indeed, for fixed degree $n$, we have 
$A_n^T(X^T)=A_n(X^T\times U/T)$ 
where $U$ is an open $T$-invariant subset of a $T$-module, such that
the quotient $U\to U/T$ exists and is a principal $T$-bundle, and that
the codimension of $V\setminus U$ is large enough. Moreover, as in the
second example of Section 1, we can find $U$ such that $U/T$ is a
product of projective spaces. It follows that 
$A_*(X^T\times U/T)=A_*(X^T)\otimes A_*(U/T)$ \cite{Fulton} Chapter 3
(note that the analogue of the K\"unneth formula does not hold for
Chow groups of products of arbitrary schemes). 

We can now state a version of the localization theorem,
due to Edidin and Graham \cite{EG3} in the more general setting of
higher equivariant Chow groups.

\begin{theorem}\label{locc}
For any $T$-scheme $X$, the $S_{\Z}$-linear map 
$$
i_{T*}:A_*^T(X^T)\to A_*^T(X)
$$
becomes an isomorphism after inverting finitely many non trivial
characters.
\end{theorem}
\begin{proof} (after \cite{Br3} 2.3) By assumption, $X$ is a finite
union of $T$-stable affine open subsets $X_i$. The ideal of each fixed
point subscheme $X_i^T$ is generated by all regular functions on $X_i$
which are eigenvectors of $T$ with a non trivial weight. We can choose
a finite set of such generators $(f_{ij})$, with respective weights
$\chi_{ij}$.

Let $Y\subset X$ be a $T$-invariant subvariety. If $Y$ is not fixed
pointwise by $T$, then one of the $f_{ij}$ defines a non zero rational
function on $Y$. Thus, we can write
$$
[Y]=\chi_{ij}^{-1}div(f_{ij})
$$ 
in $A_*^T(X)[1/\chi_{ij}]$. It follows by induction that $i_{T*}$ is
surjective after inverting the $\chi_{ij}$'s.

For injectivity, we may assume that $X$ is not fixed pointwise by $T$.
Let $Y$ be an irreducible component of $X$ not contained in $X^T$.
Choose $f_{ij}$ as before; denote by $\vert D\vert$ the union of
the support of the divisor of $f_{ij}$ in $Y$, and of the irreducible
components of $X$ which are not equal to $Y$. Then $\vert D\vert$
contains all fixed points of $X$. Denote by $j:\vert D\vert\to X$ the
inclusion. Let $U\to U/T$ be as in the definition of equivariant Chow
groups, and let $L(\chi_{ij})$ be the line bundle on $U/T$ associated
with the character $\chi_{ij}$. Denote by $p_X:X\times_T U\to U/T$ the
projection. Then we have a pseudo-divisor on $X\times_T U$
\cite{Fulton} 2.2:
$$
(p_X^*L(\chi_{ij}), \vert D\vert\times_T U,f_{ij})
$$
which defines a homogeneous map of degree -1
$$
j^*:A_*^T(X)\to A_*^T(\vert D\vert)
$$
such that composition $j^*\circ j_*$ is multiplication by $\chi_{ij}$.
Thus, the map
$$
j_*:A_*^T(\vert D\vert)\to A_*^T(X)
$$
is injective after inverting $\chi_{ij}$. We conclude by Noetherian
induction.
\end{proof}

\bigskip
If moreover $X$ is smooth, then the fixed point scheme $X^T$ 
is smooth, too. Then, as in the end of Section 1, we see that the
analogue of Theorem \ref{loc} holds in the algebraic setting. 
We refer to \cite{EG3} for applications to the Bott residue formula
\cite{Bott}.

Any projective smooth algebraic $T$-variety $X$ admits a $T$-stable
cellular decomposition, by a result of Bialynicki-Birula \cite{BB}. It
follows easily that the cycle maps $cl^T_X:A^*_T(X)\to H^*_T(X,\Z)$
and $cl_X:A^*(X)\to H^*(X,\Z)$ are isomorphisms, and that the
$S_{\Z}$-module $A^*_T(X)$ is free. The following extension of this
result is proved in \cite{Br3} 3.2. 
\begin{theorem}\label{loque}
For any projective smooth $T$-variety $X$, the $S$-module
$A^*_T(X)_{\Q}$ is free. Furthermore, the pull-back map
$$
i_T^*:A^*_T(X)_{\Q}\to A^*_T(X^T)_{\Q}
$$
is injective, and becomes surjective after inverting finitely many
non trivial characters of $T$. If moreover the cycle map
$$
cl_{X^T}:A^*(X^T)_{\Q}\to H^*(X^T,\Q)
$$
is an isomorphism, then both cycle maps
$$
cl_X^T:A^*_T(X)_{\Q}\to H^*_T(X),~cl_X:A^*(X)_{\Q}\to H^*(X)
$$
are isomorphisms, too.
\end{theorem}

It follows that the analogue of Theorem \ref{locprec} holds in this
setting, with the same proof.

\bigskip
 
Denote by $Q$ the quotient field of $S_{\Z}$; then $Q$ is the field of
rational functions in $\dim(T)$ variables, with rational coefficients.

\begin{corollary}\label{mult} Let $X$ be an equidimensional $T$-scheme
with finite fixed point set. Then we have in 
$A_*^T(X)\otimes_{S_{\Z}} Q\simeq A_*^T(X^T)\otimes_{S_{\Z}} Q
\simeq \oplus_{x\in X^T} Q[x]$:
$$
[X]=\sum_{x\in X^T} (e_x X) [x]
$$
for uniquely defined coefficients $e_x X\in Q$ which are homogeneous
of degree $-\dim(X)$. If moreover $X$ is smooth at $x$, then 
$$e_x X=\frac{1}{det(T_x X)}.$$
\end{corollary}
\begin{proof}
The first assertion follows from Theorem \ref{locc}. For the second
one, we may replace $X$ by an open $T$-stable neighborhood of $x$, and
thus assume that $X$ is smooth. Denote by $i_x$ the inclusion of
$x$ into $X$. Then we know that $i_x^*\circ i_{x*}$ is multiplication
by the determinant of the action of $T$ in $T_x X$. Comparing
coefficients of $[x]$ in $i_x^*[X]=[x]$, we obtain our result.
\end{proof}
\bigskip

The homogeneous rational function $e_x X$ is called the 
{\sl equivariant multiplicity} of $X$ at $x$. The following
observation allows to compute it in many examples.

\begin{lemma}\label{resolution}
Let $X$, $Y$ be $T$-varieties and let $\pi:Y\to X$ be a proper
surjective equivariant morphism of finite degree $d$. If $Y^T$ is
finite, then we have for any $x\in X^T$:
$$
e_x X=\frac{1}{d}\sum_{y\in Y^T,\pi(y)=x} e_y Y.
$$
\end{lemma}

This follows immediately from equality $\pi_*[Y]=d[X]$ in $A_*^T(X)$.
In particular, if $\pi$ is a resolution of singularities (that is, $Y$
is smooth and $\pi$ is birational), then
$$
e_x X=\sum_{y\in Y^T, \pi(y)=x} \frac{1}{det(T_yY)}.
$$

\noindent
{\bf Example} (Schubert varieties) 

Recall that the $T$-fixed points in the flag variety $G/B$ are the
$wB$ for $w\in W$. Moreover, the $T$-fixed points in the Schubert
variety
$$
X(w):=\overline{BwB}/B
$$
are the $xB$ where $x\in W$ and $x\leq w$ for the Bruhat ordering.
Because the map $w\mapsto wB$ is injective, we will write $w$ for
$wB$, etc. Let us determine $e_x X(w)$.

Denote by $\Phi$ the root system of $(G,T)$, by $\Phi^+$ the subset of
positive roots such that the roots of $(B,T)$ are {\sl negative}, and
by $\Pi$ the corresponding set of simple roots. Then the set of
weights of $T$ in the tangent space $T_w G/B$ is $w(\Phi^+)$ and each
such weight has multiplicity one. We thus have
$$
e_w G/B=\prod_{\alpha\in w(\Phi^+)}\alpha^{-1}.
$$
Moreover, $w$ is a smooth point of $X(w)$ and the set of weights of
$T_w X(w)$ is $\Phi^-\cap w(\Phi^+)$, whence
$$
e_w X(w)=\prod_{\alpha\in\Phi^-\cap w(\Phi^+)}\alpha^{-1}.
$$
To obtain $e_x X(w)$ for arbitrary $x$, we use the Bott-Samelson
resolution of singularities of $X(w)$, which we recall briefly. For
any $\alpha\in\Pi$, denote by $s_{\alpha}$ the corresponding
reflection and by $P_{\alpha}=B\cup Bs_{\alpha}B$ the corresponding
minimal parabolic subgroup. We can find $\alpha\in\Pi$ and $\tau\in W$
such that $w=s_{\alpha}\tau$ and that $\ell(w)=\ell(\tau)+1$ (here
$\ell$ denotes the length function on $W$). Then
$X(w)=P_{\alpha}X(\tau)$ and the map
$$
\pi:P_{\alpha}\times_B X(\tau)\to X(w)
$$
is birational. Iterating this construction, we obtain a
$B$-equivariant birational map
$$
P_{\alpha_1}\times_B P_{\alpha_2}\times_B\cdots\times_B P_{\alpha_n}/B
\to X(w)
$$
where $w=s_{\alpha_1}s_{\alpha_2}\cdots s_{\alpha_n}$ is a reduced
decomposition. Moreover, the variety on the left is projective, smooth
and contains only finitely many $T$-fixed points: the classes of
sequences $(s_1,\ldots,s_n)$ where each $s_j$ is either $s_{\alpha_j}$
or 1.

We claim that
$$
e_x X(w)=
\frac{e_x X(\tau)-s_{\alpha}(e_{s_{\alpha}x}X(\tau))}{\alpha},
$$
which determines inductively equivariant multiplicities of Schubert
varieties. Indeed, the $T$-fixed points in the fiber $\pi^{-1}(x)$ are
the classes $(1,x)B$ and $(s_{\alpha},s_{\alpha}x)B$, where the first
(resp. second) point occurs if and only if $x\leq\tau$ (resp. 
$s_{\alpha}x\leq\tau$).
Finally, $(1,x)B$ has a $T$-invariant neighborhood isomorphic to
$\C(\alpha)\times X(\tau)$ where $\C(\alpha)$ is the one-dimensional
$T$-module with weight $\alpha$. It follows that 
$$
e_{(1,x)B} P_{\alpha}\times_B X(\tau)=\frac{e_x X(\tau)}{\alpha}
$$
and that a similar equality holds for 
$e_{(s_{\alpha},s_{\alpha}x)B} P_{\alpha}\times_B X(\tau)$. Applying
Lemma \ref{resolution}, we obtain our claim.

Let now
$$
w=s_{\alpha_1}\cdots s_{\alpha_n}
$$
be a reduced decomposition. Applying repeatedly the formula above (or
equivalently, using the full Bott-Samelson resolution), we obtain an
explicit expression for $e_x X(w)$ \cite{Ar2} Proposition 3.3.1 and
\cite{Ross}:
$$
e_x X(w)=\sum_{s_1,\ldots,s_n}
\prod_{j=1}^n s_1\cdots s_j(\alpha_j^{-1})
$$
where the sum runs over all sequences $(s_1,\ldots,s_n)$ such that
$s_j=s_{\alpha_j}$ or 1, and that $s_1\cdots s_n=x$ (such a sequence
is called a {\sl subexpression} of the reduced expression for $w$).

Consider for example $w=s_{\alpha}s_{\beta}$ where $\alpha$ and
$\beta$ are distinct simple roots. Then we obtain
$$
e_1 X(s_{\alpha}s_{\beta})=
-e_{s_{\beta}} X(s_{\alpha}s_{\beta})=\frac{1}{\alpha\beta},~
e_{s_{\alpha}s_{\beta}} X(s_{\alpha}s_{\beta})=
-e_{s_{\alpha}} X(s_{\alpha}s_{\beta})=
\frac{1}{\alpha s_{\alpha}(\beta)}.
$$
If moreover $\alpha$ and $\beta$ are connected in the Dynkin diagram,
then $s_{\alpha}s_{\beta}s_{\alpha}$ has length 3, and we obtain
$$
e_1 X(s_{\alpha}s_{\beta}s_{\alpha})=
-e_{s_{\alpha}} X(s_{\alpha}s_{\beta}s_{\alpha})=
\frac{-\langle \beta,\alpha^{\vee}\rangle}
{\alpha\beta s_{\alpha}(\beta)}
$$
which shows that $e_x X(w)$ is not always the inverse of a product of
roots.

To understand better $e_x X(w)$, we recall the construction of a slice
${\cal N}_{x,w}$ to the orbit $Bx$ in $X(w)$. Denote by $U$, $U_x$
(resp. $U^-$, $U^-_x$) the unipotent subgroup of $G$ normalized by $T$
with root set $\Phi^-$, $\Phi^-\cap x(\Phi^+)$ (resp. $\Phi^+$,
$\Phi^+\cap x(\Phi^+)$). Then, by the Bruhat decomposition, the set 
$$
U^-xB/B:=U^-x
$$
is an open $T$-invariant neighborhood of $x$ in $G/B$, and the product
map
$$
U_x\times U^-_x x\to U^-x,~(u,\xi)\mapsto u\xi
$$
is an isomorphism, which restricts to an isomorphism 
$U_x\times x\to Bx$. Thus, 
$$
{\cal N}_{x,w}:=U^-_x\cap X(w)
$$
is a locally closed, $T$-stable subvariety of $X(w)$ containing $x$,
and the product map
$$
U_x\times{\cal N}_{x,w}\to X(w)
$$
is an open immersion. It follows that
$$
e_x X(w)=(e_x U_x)(e_x{\cal N}_{x,w})=
(\prod_{\alpha\in\Phi^-\cap x(\Phi^+)}\alpha^{-1})e_x{\cal N}_{x,w}.
$$
Moreover, $e_x{\cal N}_{x,w}$ is a rational function of degree
$-\dim{\cal N}_{x,w}=l(x)-l(w)$.

Consider the simplest case where $l(w)=l(x)+1$. Choose a reduced
decomposition $w=s_{\alpha_1}\cdots s_{\alpha_n}$. Then there exists
an index $j$ such that 
$$
x=s_{\alpha_1}\cdots \widehat{s_{\alpha_j}}\cdots s_{\alpha_n}.
$$
It follows that $w=s_{\beta}x$ where 
$\beta=s_{\alpha_1}\cdots s_{\alpha_{j-1}}(\alpha_j)$ is in
$\Phi^+\cap x(\Phi^+)$. Now the explicit formula for $e_x X(w)$ gives
immediately 
$$
e_x X(w)=(\prod_{\alpha\in\Phi^-\cap x(\Phi^+)}\alpha^{-1})\beta^{-1},~
e_x{\cal N}_{x,w}=\beta^{-1}.
$$
In fact, the $T$-variety ${\cal N}_{x,w}$ is isomorphic to
the module with weight $\beta$ (this will be checked at the end of
Section 5).

For arbitrary $x$ and $w$, we are led to the following
\bigskip

\noindent
{\bf Question.} What is the structure of the $T$-variety 
${\cal N}_{x,w}$?
\smallskip

This is especially interesting in the case where $X(x)$ is an
irreducible component of the singular locus of $X(w)$; equivalently,
$x$ is an isolated singular point of ${\cal N}_{x,w}$. Then the
singularity of ${\cal N}_{x,w}$ at $x$ is the generic singularity of
$X(w)$ along $Bx$. 

Little seems to be known about this question. Here are two partial
results: In the case where $l(w)=l(x)+2$, it can be shown that 
${\cal N}_{x,w}$ is a toric surface for a quotient of $T$ (see the end
of Section 5), and thus, a cone over the projective line.
On the other hand, generic singularities of Schubert varieties in
$G/P$ are described in \cite{BP} when $P\supset B$ is a maximal
parabolic subgroup associated with a minuscule or cominuscule
fundamental weight. These singularities turn out to be multicones over
homogeneous spaces.

\bigskip

Back to the general case of a torus action with isolated fixed points,
let us give an interpretation of equivariant multiplicity for an
attractive fixed point, in the following
sense.
\bigskip

\noindent
{\bf Definition.} A fixed point $x$ is {\sl attractive} if all weights
of $T$ in the tangent space $T_x X$ are contained in some open
half-space, that is, some one-parameter subgroup of $T$ acts on 
$T_x X$ with positive weights only.
\bigskip

Equivalently, all $T$-orbits in some open $T$-stable neighborhood of
$x$ contain that point in their closure. In fact, the set
$$
X_x=:\{\xi\in X~\vert~x\in\overline{T\xi}\}
$$
is the unique open affine $T$-stable neighborhood of $x$ in $X$;
clearly, $x$ is the unique closed $T$-orbit in $X_x$.

For attractive $x$, we denote by $\chi_1,\ldots,\chi_n$ the weights of
$T$ in $T_x X$. Let $\Xi_*(T)$ be the lattice of one-parameter
subgroups of $T$, and $\Xi_*(T)_{\R}$ the associated real vector space.
Set
$$
\sigma_x:=\{\lambda\in \Xi_*(T)_{\R}~\vert~
\langle\lambda,\chi_i\rangle\geq 0~{\rm for}~1\leq i\leq n\}.
$$
Then $\sigma_x$ is a rational polyhedral convex cone in $\Xi_*(T)_{\R}$
with a non empty interior $\sigma_x^0$. Any 
$\lambda\in\sigma_x^0\cap \Xi_*(T)$ defines a grading of the algebra of
regular functions $\C[X_x]$, by setting
$$
\C[X_x]_n:=
\bigoplus_{\chi\in \Xi(T),\langle\chi,\lambda\rangle=n} \C[X_x]_{\chi}.
$$
Let $d$ be the dimension of $X$. Then there exists a positive rational
number $e$ such that
$$
\sum_{m=0}^n \dim \C[X_x]_m=e\frac{n^d}{d!}+o(n^d).
$$ 
This number is called the multiplicity of $\C[X_x]$ for the grading
defined by $\lambda$. It is easily seen that $e$ is the value at
$\lambda$ of $e_x X$  (viewed as a rational function on
$\Xi_*(T)_{\R}$), and that the product $\chi_1\cdots\chi_n=det(T_x X)$
is a denominator for $e_x X$ \cite{Br3} 4.4. The homogeneous
polynomial function 
$$
\chi_1\cdots\chi_n e_x X = det(T_x X) e_x X=:J_x X
$$
of degree $n-\dim(X)$ is the {\sl Joseph polynomial} introduced in
\cite{Joseph} in relation to representation theory, see also
\cite{BBM}.
\bigskip

In the example of the flag manifold, each fixed point $x$ is
attractive, and $(G/B)_x=xU^-B/B$ where $U^-\subset G$ is the
unipotent subgroup normalized by $T$ with 
root set $\Phi^+$. Moreover, the cone $\sigma_x$ is the image of the
positive Weyl chamber under $x^{-1}\in W$. Because ${\cal N}_{x,w}$ is
contained in $U^-_x$, the product of all roots in 
$\Phi^+\cap x(\Phi^+)$ is a denominator for $e_x{\cal N}_{x,w}$. A
more precise result will be given at the end of Section 5.

On the other hand, if $X$ is a toric variety, then any fixed point
$x\in X$ is attractive, and $\sigma_x$ is the cone associated with the
affine toric variety $X_x$. Moreover, for any
$\lambda\in \sigma_x^0$, the value at $\lambda$ of $e_x X$ is $d!$
times the volume of the convex polytope 
$P_x(\lambda):=\{x\in\Xi(T)_{\R}~\vert~
\langle\mu,x\rangle\geq 0~\forall\mu\in\sigma_x,~
\langle\lambda,x\rangle\leq 1\}
$ (see \cite{Br3} 5.2). Here the volume form on $\Xi(T)_{\R}$ is
normalized so that the quotient by the lattice $\Xi(T)$ has volume 1.

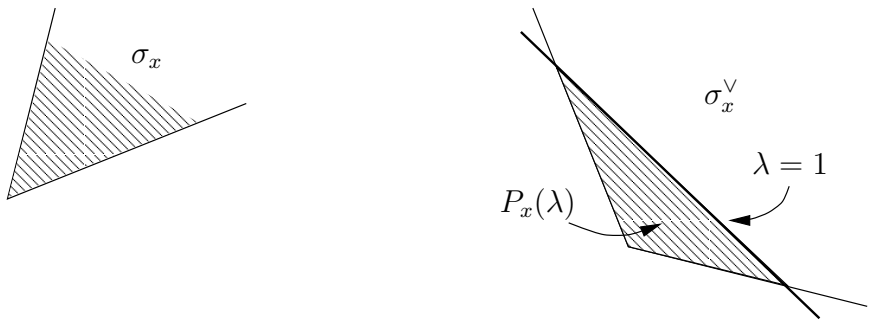
\begin{center}
\begin{picture}(0,0)%
\epsfig{file=figure4.pstex}%
\end{picture}%
\setlength{\unitlength}{0.00083300in}%
\begingroup\makeatletter\ifx\SetFigFont\undefined%
\gdef\SetFigFont#1#2#3#4#5{%
  \reset@font\fontsize{#1}{#2pt}%
  \fontfamily{#3}\fontseries{#4}\fontshape{#5}%
  \selectfont}%
\fi\endgroup%
\begin{picture}(5424,1984)(1189,-1433)
\put(5551,-61){\makebox(0,0)[lb]{\smash{\SetFigFont{12}{14.4}{\rmdefault}{\mddefault}{\updefault}$\sigma_x^\vee$}}}
\put(4276,-736){\makebox(0,0)[lb]{\smash{\SetFigFont{12}{14.4}{\rmdefault}{\mddefault}{\updefault}$P_x(\lambda)$}}}
\put(5851,-501){\makebox(0,0)[lb]{\smash{\SetFigFont{12}{14.4}{\rmdefault}{\mddefault}{\updefault}$\lambda=1$}}}
\put(1976,194){\makebox(0,0)[lb]{\smash{\SetFigFont{12}{14.4}{\rmdefault}{\mddefault}{\updefault}$\sigma_x$}}}
\end{picture}

\end{center}

\section{Criteria for (rational) smoothness}

We begin by giving a smoothness criterion at an attractive fixed point
of a torus action, in terms of fixed points of subtori of codimension
1 as in Theorem \ref{locprec}.

\begin{theorem}\label{sm}
For a $T$-scheme $X$ with an attractive fixed point $x$, the following
conditions are equivalent:

(i) $X$ is smooth at $x$.

(ii) For any subtorus $T'\subset T$ of codimension 1, the fixed point
set $X^{T'}$ is smooth at $x$, and we have
$$
e_x X=\prod_{T'} e_x(X^{T'})
$$
(product over all subtori of codimension 1).
\end{theorem}
\begin{proof}
Replacing $X$ by $X_x$, we may assume that $X$ is affine.

(i)$\Rightarrow$(ii) By the graded Nakayama lemma, $X$ is
equivariantly isomorphic to a $T$-module $V$. Then each $X^{T'}$ is
isomorphic to the $T$-module $V^{T'}$; thus, $X^{T'}$ is smooth.
Moreover, we have by Corollary \ref{mult}:
$$
e_x X=\prod\frac{1}{\chi^{\dim(V_{\chi})}}
$$
(product over all characters $\chi$ of $T$, where $V_{\chi}$ denotes
the corresponding eigenspace). Then
$$
e_x(X^{T'})=\prod_{\chi,\chi\vert T'=0}\frac{1}{\chi^{\dim(V_{\chi})}}
$$
for each subtorus $T'$. Our formula follows.

(ii)$\Rightarrow$(i) There exists an equivariant closed embedding
$$
\iota:X\to T_x X
$$
such that $\iota(x)=0$. For any subtorus $T'\subset T$ of codimension
one, we have a subspace $T_x(X^{T'})$ of $T_x X$. Let $V\subset T_x X$
be the span of all these subspaces; then $V$ is a $T$-submodule of
$T_x X$. Choose an equivariant projection 
$$
p:T_x X\to V
$$
and denote by 
$$
\pi:X\to V
$$
the composition of $\iota$ and $p$. Then $\pi$ is equivariant, $\pi(x)=0$
and the differential of $\pi$ at $x$ induces isomorphisms
$T_x(X^{T'})\to V^{T'}$ for all subtori $T'\subset T$ of codimension 1.
Because $X^{T'}$ is smooth at the attractive point $x$, it follows
by the graded Nakayama lemma that $\pi$ restricts to isomorphisms
$X^{T'}\to V^{T'}$. 

We claim that the morphism $\pi$ is finite. By the graded Nakayama
lemma again, it is enough to check that the set $\pi^{-1}(0)$ consists
of $x$, or even that $\pi^{-1}(0)$ is finite (because $x$ is the
unique fixed point of $X$). But if $\pi^{-1}(0)$ is infinite, then
this closed $T$-stable subset of $X$ contains a $T$-stable closed
curve $C$. So $C$ contains $x$ and is fixed pointwise by a subtorus
$T'\subset T$ of codimension 1. This contradicts the fact that
restriction of $\pi$ to $X^{T'}$ is injective.

Now observe that
$$
\dim(X)=-\deg(e_x X)=-\sum_{T'}\deg(e_x(X^{T'}))=\sum_{T'}\dim(X^{T'})
=\dim(V).
$$
Thus, the finite morphism $\pi:X\to V$ is surjective. Then the algebra
of regular functions $\C[X]$ is a finite torsion-free module over
$\C[V]$; let $d$ be its rank. We have
$$
e_x X=d e_x V=d\prod_{T'}e_x(V^{T'})=d\prod_{T'} e_x(X^{T'})=d e_x X
$$
and $e_x X$ is non zero, because no $e_x(X^{T'})$ is. Thus, $d=1$, that
is, $\pi$ is birational. As $V$ is normal, it follows that $\pi$ is an
isomorphism.
\end{proof}
\bigskip

We now adapt these arguments to obtain a similar criterion for 
rational smoothness at an attractive fixed point. Recall the
following
\bigskip

\noindent
{\bf Definition} \cite{KL} An algebraic variety $X$ is {\sl rationally
smooth of dimension $n$} (or a rational cohomology manifold of
dimension $n$) if we have $H^m(X,X\setminus x)=0$ for $m\neq n$, and
$H^n(X,X\setminus x)\simeq\Q$, for all $x\in X$. A point $x\in X$ is
rationally smooth if it admits a rationally smooth neighborhood.
\bigskip

By recent work of A. Arabia \cite{Ar3}, an attractive fixed point $x$
such that all weights in $T_x X$ have multiplicity 1 is rationally
smooth if and only if: a punctured neighborhood of $x$ is rationally
smooth, and $e_x X$ is the inverse of a polynomial. Here is an
extension of this result. 

\begin{theorem}\label{rsm}
Let $X$ be a $T$-variety with an attractive fixed point $x$ such
that a punctured neighborhood of $x$ in $X$ is rationally smooth. Then
the following conditions are equivalent:

(i) The point $x$ is rationally smooth.

(ii) For any subtorus $T'\subset T$ of codimension 1, the point $x$ is
rationally smooth in $X^{T'}$, and there exists a positive
rational number $c$ such that
$$
e_x X=c\prod_{T'} e_x(X^{T'})
$$
(product over all subtori of codimension 1). If moreover each
$X^{T'}$ is smooth, then $c$ is an integer.

(iii) For any subtorus $T'\subset T$ of codimension 1, the point $x$
is rationally smooth in $X^{T'}$, and we have
$\dim(X)=\sum_{T'}\dim(X^{T'})$ (sum over all subtori of codimension
1).
\end{theorem}
\begin{proof}
We begin with some observations on the local structure of $X$ at $x$.
First, we may assume that $X$ is affine. Then $x$ is the unique closed
$T$-orbit in $X$, and the punctured space $\dot X:=X\setminus x$ is
rationally smooth. 

Choose a one-parameter subgroup $\lambda$ such that all weights of
$\lambda$ in $T_x X$ are positive. For the action of $\C^*$ on
$\dot X$ via $\lambda$, the quotient
$$
\dot X/\C^*=:\CP(X)
$$
exists and is a projective variety. Indeed, let $\iota:X\to T_x X$ be
an equivariant embedding such that $\iota(x)=0$. Then $\iota$ induces
a closed embedding of $\CP(X)$ into $\CP(T_x X)$, and $\CP(T_x X)$ is a
weighted projective space. 

Moreover, $\CP(X)$ is rationally smooth. Indeed, $\dot X$ is
covered by $\C^*$-stable open subsets of the form
$$\C^*\times_{\Gamma} Y
$$ 
where $\Gamma$ is a finite subgroup of $\C^*$, and $Y$ is a
$\Gamma$-stable subvariety of $\dot X$; then $\CP(X)$ is covered by
the quotients $Y/\Gamma$. Because $\dot X$ is rationally smooth,
$\C^*\times Y$ and $Y$ are, too, and so are the quotients $Y/\Gamma$.

The action of $T$ on $X$ induces an action of $T/\lambda(\C^*)$ on
$\CP(X)$ for which the fixed point set is the disjoint union of the
$\CP(X^{T'})$. Indeed, $T$-fixed points in $\CP(X)$ correspond to
$T$-orbits of dimension 1 in $X$.

Observe that $x$ is rationally smooth of dimension $n$ if and only if
$\dot X$ is a rational cohomology sphere of dimension $n-1$.
Indeed, because the action of $\C^*$ on $X$ extends to a map 
$\C\times X\to X$ sending $0\times X$ to 0, the space $X$ is
contractible. Thus, $H^m(X)=0$ for all $m>0$. Now our claim follows
from the long exact sequence in relative cohomology.

Observe finally that $x$ is rationally smooth if and
only if $\CP(X)$ is a rational cohomology complex
projective space $\CP^{d-1}$ where $d=\dim(X)$. Indeed, the group
$S^1\subset \C^*$ acts on $\dot X$ without fixed points,
whence a Gysin exact sequence
$$
\cdots\to H^m(\dot X)\to H^{m-1}(\dot X/S^1)
\to H^{m+1}(\dot X/S^1)\to H^{m+1}(\dot X)
\to\cdots
$$
It follows that $\dot X$ is a rational cohomology $(n-1)$-sphere
if and only if: $n$ is even and $\dot X/S^1$ has the rational
cohomology of $\CP^{(n-1)/2}$. But the map 
$$
\dot X/S^1\to\dot X/\C^*=\CP(X)
$$ 
induces an isomorphism in rational cohomology. Moreover, because 
$\CP(X)$ is rationally smooth and projective of dimension $d-1$, 
the maximal degree occuring in its rational cohomology is $2(d-1)$
and this forces $n=2d$.

(i)$\Rightarrow$(ii) The first assertion follows from the discussion
above, together with a theorem of Smith: the fixed point set of a
torus acting on a rational cohomology sphere is a rational cohomology
sphere as well \cite{Hsiang} Theorem IV.2.

For each $T'$, choose a homogeneous system of parameters of the
algebra $\C[X^{T'}]$ (graded by the action of $T/T'\simeq\C^*$). This
defines a $T$-equivariant finite surjective morphism 
$$
\pi_{T'}:X^{T'}\to V^{T'}
$$
where $V^{T'}$ is a $T$-module with trivial action of $T'$. We can
extend $\pi_{T'}$ to an equivariant morphism 
$$
X\to V^{T'}
$$
still denoted by $\pi_{T'}$. Thus, we obtain an equivariant morphism
$$
\pi:X\to V
$$
where $V$ is the direct sum of the $V^{T'}$ over all $T'$ (observe
that only finitely many such spaces are non zero). Moreover, the
morphism $\pi$ is finite. Indeed, it restricts to a finite morphism on
each $X^{T'}$, and thus $\pi^{-1}(0)$ contains no $T$-stable curve by
the argument of Theorem \ref{sm}.

Because $\CP(X)$ and all $\CP(X^{T'})$ are rational cohomology
projective spaces, we have
$$\displaylines{
\dim(X)=\chi(\CP(X))=\chi(\CP(X)^T)
=\sum_{T'}\chi(\CP(X^{T'}))
\hfill\cr\hfill
=\sum_{T'}\dim(X^{T'})=\sum_{T'}\dim(V^{T'})=\dim(V).
\cr}$$
Thus, $\pi$ is surjective. Let $d$ be its degree, and let $d_{T'}$ be
the degree of $\pi_{T'}$. Then we have by Lemma \ref{resolution}
$$
e_x X=d e_0 V=d\prod_{T'}e_0(V^{T'})=
\frac{d}{\prod_{T'} d_{T'}}\prod_{T'}e_x(X^{T'}).
$$
If moreover each $X^{T'}$ is smooth, then we can take each $\pi_{T'}$
to be the identity. So each $d_{T'}$ is 1, and $c=d$ is an integer.
\bigskip

(ii)$\Rightarrow$(iii) Because $e_x X$ is homogeneous of degree
$-\dim(X)$, we obtain $\dim(X)=\sum_{T'}\dim(X^{T'})$.
\bigskip

(iii)$\Rightarrow$(i) Because $\CP(X)$ is projective and rationally
smooth, the spectral sequence associated with the fibration 
$\CP(X)\times_T E_T\to B_T$ degenerates (by the criterion of Deligne,
see e.g. \cite{Joshua}). Thus, the $S$-module $H^*_T(\CP(X))$ is free,
and the map 
$$
\rho:H^*_T(\CP(X))/(S_+)\to H^*(\CP(X))
$$ 
is an isomorphism. 

On the other hand, $\CP(X)^T$ has no odd cohomology
because it is the disjoint union of the $\CP(X^{T'})$, and each
$X^{T'}$ is rationally smooth. Thus, $\CP(X)^T$ has no odd equivariant
cohomology as well. By the localization theorem, the odd equivariant
cohomology of $\CP(X)$ vanishes; thus, the cohomology of $\CP(X)$ is
concentrated in even degree.

Moreover, we have
$$\displaylines{
\chi(\CP(X))={\rm rk}_S H^*_T(\CP(X))={\rm rk}_S H^*_T(\CP(X)^T)=
\chi(\CP(X)^T)
\hfill\cr\hfill
=\sum_{T'}\chi(\CP(X^{T'}))=\sum_{T'}\dim(X^{T'})=\dim(X).
\cr}$$
Thus, $\CP(X)$ has the rational cohomology of projective space of
dimension $\dim(X)-1$.
\end{proof}
\bigskip

\noindent
{\bf Example} (Schubert varieties, continued) 

Let $w,x$ in $W$ such that $x\leq w$. Notation being as in the
previous section, rational smoothness of $X(w)$ at $x$ is equivalent
to $P_{y,w}=1$ for all  $y\in W$ such that $x\leq y\leq w$, where
$P_{y,w}$ denotes the Kazhdan-Lusztig polynomial \cite{KL} Theorem A2
(this can be checked as in the beginning of the proof of Theorem
\ref{rsm}).

Set 
$$
\Phi(x,w):=\{\alpha\in x(\Phi^+)~\vert~s_{\alpha}x\leq w\}.
$$
Observe that $\Phi(x,w)$ contains $\Phi^-\cap x(\Phi^+)$
(which is the set of all $\alpha\in x(\Phi^+)$ such that
$s_{\alpha}x\leq x$), with complement the set of all 
$\alpha\in x(\Phi^+)$ such that $x<s_{\alpha}x\leq w$. Indeed, for any
$\alpha\in\Phi$, we have either $s_{\alpha}x < x$ or 
$x < s_{\alpha}x$. 

We can now state the following result, where (i) is due to A. Arabia
\cite{Ar3}, (ii) and (iii) to S. Kumar \cite{Kumar}, and (iv) to J.
Carrell and D. Peterson \cite{C}.
\begin{corollary}

\noindent
(i) There exists a homogeneous polynomial $J(x,w)\in S_{\Z}$ such that
$$
e_x X(w)=J(x,w)\prod_{\alpha\in\Phi(x,w)}\alpha^{-1}.
$$

\noindent
(ii) The point $x$ is smooth in $X(w)$ if and only if $J(x,w)=1$, that
is,
$$
e_x X(w)=\prod_{\alpha\in\Phi(x,w)}\alpha^{-1}.
$$

\noindent
(iii) The point $x$ is rationally smooth in $X(w)$ if and only if: For
any $y\in W$ such that $x\leq y < w$, the polynomial $J(y,w)$ is
constant, that is, there exists a positive integer $d(y,w)$ such that 
$$
e_y X(w)=d(y,w)\prod_{\alpha\in\Phi(y,w)}\alpha^{-1}.
$$

\noindent
(iv) The set of weights of $T$ in the tangent space $T_x X(w)$
contains $\Phi(x,w)$, and we have 
$$
l(w)\leq\vert\Phi(x,w)\vert\leq\dim T_x X(w).
$$
Moreover, $x$ is rationally smooth in $X(w)$ if and only if: For any
$y\in W$ such that $x\leq y < w$, we have $l(w)=\vert\Phi(y,w)\vert$.
Finally, $x$ is smooth in $X(w)$ if and only if: $x$ is rationally
smooth, and $\vert\Phi(x,w)\vert=\dim T_x X(w)$.
\end{corollary}
\begin{proof} We begin by proving (ii).
Let $T'\subset T$ be a subtorus of codimension 1. Then
$X(w)_x^{T'}\neq x$ if and only if $T'=ker(\alpha)$ for some
$\alpha\in\Phi(x,w)$; in this case, $X(w)_x^{T'}$ is a smooth curve,
isomorphic to the one-dimensional $T$-module with weight $\alpha$ (see
the notes of T. A. Springer in this volume, or Example 1 in Section
2). In particular, $e_x X(w)^{ker(\alpha)}=\alpha^{-1}$. 
Thus, (ii) follows from Theorem \ref{sm}. 

Now we prove (i). By the argument of Theorem \ref{sm}, there exists a
finite $T$-equivariant morphism 
$$
\pi:X(w)_x\to V
$$
where $V$ is a $T$-module with weight set $\Phi(x,w)$ (for a concrete
construction of $\pi$, embed $X(w)_x$ into $xU^-x^{-1}$; the latter is
$T$-equivariantly isomorphic to a module with weight set $x(\Phi^+)$,
which projects onto $V$). Let $d$ be the degree of $\pi$, and let
$Y\subset V$ be its image. Then we have $\pi_*[X(w)_x]=d[Y]$ in
$A_*^T(V)$, whence $e_x X(w)=de_0 Y$. Moreover, because $V$ is a
$T$-module, the $S_{\Z}$-module $A_*^T(V)$ is freely generated by
$[V]$. Write $[Y]=J(Y)[V]$ with $J(Y)\in S_{\Z}$, then 
$$
e_0 Y=J(Y)e_0 V=J(Y)\prod_{\alpha\in\Phi(x,w)}\alpha^{-1}
$$
whence (i) with $J(x,w)=dJ(Y)$.

The implication ($\Leftarrow$) in (iii) follows from Theorem 
\ref{rsm}. For ($\Rightarrow$), we use the slice ${\cal N}_{x,w}$ to
the orbit $Bx$ in $X(w)$ whose construction was recalled in the
previous section. The product map
$$
U_x\times {\cal N}_{x,w}\to X(w)
$$
is an open immersion with image $X(w)_x$.
Because $U_x$ is an affine space, $x$ is rationally smooth in $X(w)$
if and only if it is in ${\cal N}_{x,w}$. Now we argue by induction
over $l(w)-l(x)=\dim{\cal N}_{x,w}$. If this number is zero, then
$x=w$ and $X(w)$ is smooth at $x$. In the general case, the induction
hypothesis tells us that $X(w)$ is rationally smooth at $y$ for all
$y\in W$ such that $x < y\leq w$. Because the set of rationally smooth
points of $X(w)$ is open and $B$-stable, it follows that 
$X(w)\setminus\overline{Bx}$ is rationally smooth. Thus, 
${\cal N}_{x,w}\setminus x$ is rationally smooth as well, and Theorem
\ref{rsm} applies.

We already saw that $\Phi(x,w)$ is the set of weights of tangent
spaces to $x$ at $T$-stable curves in $X(w)$. This proves the first
assertion of (iv), and the inequality 
$\vert\Phi(x,w)\vert\leq \dim T_x X(w)$ as well.  

The inequality $l(w)\leq\vert\Phi(x,w)\vert$ follows from finiteness
of the morphism $\pi:X(w)_x\to V$. Moreover, equality holds if and
only if $\pi$ is surjective, or equivalently,
$$
e_x X(w)=d\prod_{\alpha\in \Phi(x,w)}\alpha^{-1}
$$
where $d$ is the degree of $\pi$. So the criterion for rational
smoothness follows from (iii). 

If $x$ is smooth in $X(w)$, then it is certainly rationally smooth,
and the weights of $T$ in $T_x X(w)$ are the weights of $T$-stable
curves through $x$ in $X(w)$, that is, the elements of $\Phi(x,w)$.
For the converse, observe that 
$$
\dim X(w)=\ell(w)=\vert\Phi(x,w)\vert=\dim T_x X(w).
$$
\end{proof}
\bigskip

Let us give some applications of these criteria. Consider first the
case where $l(w)=l(x)+1$. As in Section 4, write $w=s_{\beta}x$ where
$\beta\in \Phi^+\cap x(\Phi^+)$. Then $\Phi(x,w)$ consists of 
$\Phi^-\cap x(\Phi^+)$ and of $\beta$; moreover, 
$e_x{\cal N}_{x,w}=\beta^{-1}$. By (ii), it follows that $x$ is a
smooth point of $X(w)$, that is, the $T$-variety ${\cal N}_{x,w}$ is
isomorphic to the module with weight $\beta$. In other words, 
{\sl Schubert varieties are smooth in codimension 1}
(in fact, they are known to be normal \cite{MS}).

Consider now the case where $l(w)=l(x)+2$. Then the open Bruhat
interval $(x,w)$ contains exactly two elements $y_1,y_2$ (this follows
from the fact that, in any closed Bruhat interval, the number of
elements of even length equals the number of elements of odd length
\cite{Deo1} Proposition 3.6). Writing $y_j=s_{\beta_j}w$ as above, we
see that $\Phi(x,w)$ consists of $\Phi^-\cap x(\Phi^+)$ together with
$\beta_1$ and $\beta_2$. By (iv), it follows that 
{\sl Schubert varieties are rationally smooth in codimension 2,}
a result of Kazhdan and Lusztig \cite{KL}. 

It is easy to construct examples of Schubert varieties which are
rationally smooth, but singular in codimension 2. Indeed, let
$\alpha,\beta$ be distinct simple roots which are connected in the
Dynkin diagram. Then we saw in Section 4 that
$$
e_1 X(s_{\alpha}s_{\beta}s_{\alpha})=
-e_{s_{\alpha}} X(s_{\alpha}s_{\beta}s_{\alpha})=
\frac{-\langle \beta,\alpha^{\vee}\rangle}
{\alpha\beta s_{\alpha}(\beta)}.
$$
So the 3 dimensional Schubert variety
$X(s_{\alpha}s_{\beta}s_{\alpha})$ is rationally smooth at the
$B$-fixed point, and hence everywhere; but it is singular along the
curve $X(s_{\alpha})$ whenever 
$\langle\beta,\alpha^{\vee}\rangle\leq -2$. 

In other words, for any non simply laced group $G$, there exist
singular, rationally smooth Schubert varieties. On the other hand,
by an unpublished result of D. Peterson, rational smoothness implies
smoothness for Schubert varieties of simply laced $G$ \cite{C}. In the
case where $G={\rm SL}_n$, this is due to Deodhar \cite{Deo2} and
follows from (iv) together with the fact that the equality
$\vert\Phi(x,w)\vert=\dim T_x X(w)$ holds there for all $x$ and $w$
\cite{LS} Theorem 1.

Further applications of (rational) smoothness criteria can be found in
\cite{C} and \cite{Kumar}.


\begin{thebibliography}{100}

\bibitem{AB} M. F. Atiyah and R. Bott: {\sl The moment map and
equivariant cohomology}, Topology {\bf 23} (1984), 1-28.

\bibitem{Ar1} A. Arabia: {\sl Cycles de Schubert et cohomologie
\'equivariante de $K/T$}, Invent. math. {\bf 85} (1986), 39-52.

\bibitem{Ar2} A. Arabia: {\sl Cohomologie $T$-\'equivariante de la
vari\'et\'e des drapeaux d'un groupe de Kac-Moody}, Bull. Soc. math.
France {\bf 117} (1989), 129-165.

\bibitem{Ar3} A. Arabia: {\sl Classe d'Euler \'equivariante et points
rationnellement lisses}, preprint, December 1997.

\bibitem{Au} M. Audin: {\sl The topology of torus actions on
symplectic manifolds}, Birkh\"auser, Boston 1991.

\bibitem{BB} A. Bialynicki-Birula: {\sl Some theorems on actions of
algebraic groups}, Ann. of Math. {\bf 98} (1973), 480-497.

\bibitem{BDP} E. Bifet, C. De Concini and C. Procesi: {\sl Cohomology
of regular embeddings}, Adv. Math. {\bf 82} (1990), 1-34.

\bibitem{Borel} A. Borel: {\sl Sur la cohomologie des espaces fibr\'es
principaux et des espaces homog\`enes des groupes de Lie compacts},
Ann. of Math. {\bf 57} (1953), 115-207.

\bibitem{BBM} W. Borho, J-L. Brylinski and R. MacPherson: 
{\sl Nilpotent Orbits, Primitive Ideals, and Characteristic Classes},
Birkh\"auser, Boston 1989.

\bibitem{Bott} R. Bott: {\sl A residue formula for holomorphic vector
fields}, J. Differential Geom. {\bf 1} (1967), 311-330.

\bibitem{Bour} N. Bourbaki: {\sl Groupes et alg\`ebres de Lie},
Chapitres 7 et 8, CCLS, Paris 1975.

\bibitem{Br1} M. Brion: {\sl Sur l'image de l'application moment}, 
p. 177-192 in: S\'eminaire d'alg\`ebre Paul Dubreil et Marie-Paul
Malliavin, Lecture Notes in Math. {\bf 1296}, Springer, New York 1987.

\bibitem{Br2} M. Brion: {\sl Piecewise polynomial functions, convex
polytopes, and enumerative geometry}, in: Parameter Spaces 
(P. Pragacz, ed.), Banach Center Publications, 1996, 25-44.

\bibitem{Br3} M. Brion: {\sl Equivariant Chow groups for torus
actions}, Journal of Transformation Groups {\bf 2} (1997), 225-267.

\bibitem{BP} M. Brion and P. Polo: {\sl Generic singularities of
certain Schubert varieties}, preprint 1997.

\bibitem{C} J. Carrell: {\sl The Bruhat graph of a Coxeter group, a
conjecture of Deodhar, and rational smoothness of Schubert varieties},
pp. 53-62 in: Algebraic groups and their generalizations (W. Haboush
and B. Parshall, eds.), American Mathematical Society, Providence
1994.

\bibitem{CS} T. Chang and T. Skjelbred: {\sl The topological Schur
lemma and related results}, Ann. of Math. {\bf 100} (1974), 307-321.

\bibitem{DP} C. De Concini and C. Procesi: {\sl Complete symmetric
varieties}, pp. 1-44 in: Invariant Theory, Lecture Note in Math. 
{\bf 996}, Springer-Verlag, New York 1983.

\bibitem{Deo1} V. Deodhar: {\sl Some characterizations of the Bruhat
ordering on a Coxeter group and determination of the relative M\"obius
function}, Invent. math. {\bf 39} (1977), 187-198.

\bibitem{Deo2} V. Deodhar: {\sl Local Poincar\'e duality and
nonsingularity of Schubert varieties}, Comm. Alg. {\bf 13} (1985),
1379-1388.

\bibitem{EG1} D. Edidin and W. Graham: {\sl Characteristic classes
in the Chow ring}, J. Alg. Geom. {\bf 6} (1997), 431-443.
 
\bibitem{EG2} D. Edidin and W. Graham: {\sl Equivariant intersection
theory}, Invent. math., to appear.

\bibitem{EG3} D. Edidin and W. Graham: {\sl Localization in
equivariant intersection theory and Bott's residue formula}, Amer. J.
Math., to appear.

\bibitem{Fulton} W. Fulton: {\sl Intersection Theory}, Springer, New
York 1984.

\bibitem{FMSS} W. Fulton, R. MacPherson, F. Sottile and B. Sturmfels:
{\sl Intersection theory on spherical varieties}, J. Alg. Geom. 
{\bf 4} (1995), 181-193.

\bibitem{Gil} H. Gillet: {\sl Riemann-Roch theorems for higher
algebraic $K$-theory}, Adv. Math. {\bf 40} (1981), 203-289.

\bibitem{GKM} M. Goresky, R. Kottwitz and R. MacPherson: 
{\sl Equivariant cohomology, Koszul duality, and the localization
theorem}, Invent. math. {\bf 131} (1998), 25-84.

\bibitem{Gro} A. Grothendieck: {\sl Torsion homologique et sections
rationnelles}, dans: S\'eminaire Chevalley, Anneaux de Chow et
applications, \'Ecole Normale Sup\'erieure 1958.

\bibitem{GS} V. Guillemin and S. Sternberg: {\sl Multiplicity-free
spaces}, J. Differ. Geom. {\bf 19} (1984), 31-56.

\bibitem{Hopf} H. Hopf: {\sl \"Uber die Topologie der
Gruppen-Mannigfaltigkeiten und ihrer Verallgemeinerungen}, Ann. of
Math. {\bf 42} (1941), 22-52.

\bibitem{Hsiang} W. Y. Hsiang: {\sl Cohomology Theory of Topological
Transformation Groups}, Springer, New York 1975.

\bibitem{Hu} D. Husemoller: {\sl Fiber Bundles}, third edition,
Springer-Verlag, New York 1994.

\bibitem{I} P. Iglesias: {\sl Les SO(3)-vari\'et\'es symplectiques et
leur classification en dimension 4}, Bull. Soc. math. France {\bf 119}
(1991), 371-396.

\bibitem{Joseph} A. Joseph: {\sl On the variety of a highest weight
module}, J. Algebra {\bf 88} (1984), 238-278.

\bibitem{Joshua} R. Joshua: {\sl Vanishing of odd-dimensional
intersection cohomology}, Math. Z. {\bf 195} (1987), 239-253.

\bibitem{Kawa} K. Kawakubo: {\sl The theory of transformation groups},
Oxford University Press, Oxford 1991.

\bibitem{KL} D. Kazhdan and G. Lusztig: {\sl Representations of
Coxeter groups and Hecke algebras}, Invent. math. {\bf 53} (1979),
165-184.

\bibitem{Kirwan} F. Kirwan: {\sl Cohomology of quotients in symplectic
and algebraic geometry}, Mathematical Notes {\bf 31}, Princeton
University Press, Princeton 1984.

\bibitem{KK} B. Kostant and S. Kumar: {\sl The nil Hecke ring and
cohomology of $G/P$ for a Kac-Moody group $G$}, Adv. Math. {\bf 62}
(1986), 187-237.

\bibitem{Kumar} S. Kumar: {\sl The nil Hecke ring and singularity of
Schubert varieties}, Invent. math. {\bf 123} (1996), 471-506.

\bibitem{LS} V. Lakshmibai and C. S. Seshadri: {\sl Singular locus 
of a Schubert variety}, Bull. Amer. Math. Soc. (N.S.) {\bf 11} (1984),
363-366.

\bibitem{Lerman} E. Lerman: {\sl Symplectic cuts}, Math. Res. Lett.
{\bf 2} (1995), 247-258.

\bibitem{LP} P. Littelmann and C. Procesi: {\sl Equivariant Cohomology
of Wonderful Compactifications}, p. 219-262 in: Operator Algebras,
Unitary Representations, Enveloping Algebras, and Invariant Theory,
Birkh\"auser, Basel 1990.

\bibitem{MS} V. Mehta and V Srinivas: {\sl Normality of Schubert
varieties}, Amer. J. of Math. {\bf 109} (1987), 987-989.

\bibitem{Ross} W. Rossmann: {\sl Equivariant multiplicities on complex
varieties}, pp. 313-330 in: Orbites unipotentes et faisceaux pervers,
Ast\'erisque {\bf 173-174}, Soc. math. France 1989.

\bibitem{S} H. Sumihiro: {\sl Equivariant completion}, J. Math. Kyoto
Univ. {\bf 14} (1974), 1-28.

\bibitem{W} C. Woodward: {\sl The classification of transversal
multiplicity-free actions}, Ann. Glob. An. Geom. {\bf 14} (1996),
3-42.
\end{thebibliography}
\end{document}